\documentclass[a4paper,12pt]{amsart}

\usepackage{amsthm}
\usepackage{amsmath}
\usepackage{amsfonts}
\usepackage{amssymb}
\usepackage{latexsym}
\usepackage{times}
\usepackage{graphicx}
\usepackage{yhmath}

\theoremstyle{definition}
\newtheorem{defn}{\indent\bf Definition}
\newtheorem{rem}[defn]{\indent\bf Remark}

\theoremstyle{plain}
\newtheorem{lemma}[defn]{\indent\bf Lemma}
\newtheorem{prop}[defn]{\indent\bf Proposition}
\newtheorem{thm}[defn]{\indent\bf Theorem}
\newtheorem{cor}[defn]{\indent\bf Corollary}

\begin{document}

\title[Homogeneous quotients by a discrete group]{On the countable measure-preserving relation\\ induced on a homogeneous quotient by the action of a discrete group}
\author[Florin R\u adulescu]{Florin R\u adulescu ${}^*$
  \\ \\
Dipartimento di Matematica\\ Universit\`a degli Studi di Roma ``Tor Vergata''}

%\maketitle

\thispagestyle{empty}

\renewcommand{\thefootnote}{}
\footnotetext{${}^*$ Member of the Institute of  Mathematics ``S. Stoilow" of the Romanian Academy}
\footnotetext{${}^*$
Supported in part by PRIN-MIUR   and  by a grant of the Romanian National Authority for Scientific Research,  project number   PN-II-ID-PCE-2012-4-0201}

\footnotetext{  email radulesc@mat.uniroma2.it}
\footnotetext{
MSC classes:	11F72, 	11F25, 46L65, 	46L54}

\def\tilde{\widetilde}
\def\a{\alpha}
\def\T{\theta}
\def\PSL{\mathop{\rm PSL}\nolimits}
\def\SL{\mathop{\rm SL}\nolimits}
\def\PGL{\mathop{\rm PGL}}
\def\Per{\mathop{\rm Per}}
\def\GL{\mathop{\rm GL}}
\def\Out{\mathop{\rm Out}}
\def\Int{\mathop{\rm Int}}
\def\Aut{\mathop{\rm Aut}}
\def\ind{\mathop{\rm ind}}
\def\card{\mathop{\rm card}}
\def\d{{\rm d}}
\def\Z{\mathbb Z}
\def\R{\mathbb R}
\def\cR{{\mathcal R}}
\def\tPsi{{\tilde{\Psi}}}
\def\Q{\mathbb Q}
\def\N{\mathbb N}
\def\C{\mathbb C}
\def\bH{\mathbb H}
\def\Y{{\mathcal Y}}
\def\L{{\mathcal L}}
\def\G{{\mathcal G}}
\def\U{{\mathcal U}}
\def\H{{\mathcal H}}
\def\I{{\mathcal I}}
\def\A{{\mathcal A}}
\def\S{{\mathcal S}}
\def\O{{\mathcal O}}
\def\V{{\mathcal V}}
\def\D{{\mathcal D}}
\def\B{{\mathcal B}}
\def\K{{\mathcal K}}
\def\cC{{\mathcal C}}
\def\cR{{\mathcal R}}
\def\cX{{\mathcal X}}
\def\ptimes{\mathop{\boxtimes}\limits}
\def\potimes{\mathop{\otimes}\limits}
\def\piK{\pi_{\operatorname{Koop}}}

\begin{abstract}
We consider a countable discrete group $G$
acting ergodicaly and a.e. freely, by measure-preserving transformations, on
an infinite measure space $(\mathcal X,\nu)$ with $\sigma$-finite
measure $\nu$. Let  $\Gamma \subseteq G$ be an almost
normal subgroup with fundamental domain $F\subseteq \cX $ of finite measure.  Let $\mathcal R_G$ be the  countable measurable equivalence relation
 on $\mathcal X$  determined by the orbits of $G$. Let $\mathcal R_G| F$
be its restriction to $F$.

We find an explicit presentation, by  generators and relations, for the von Neumann algebra associated, by the Feldman-Moore (\cite{FM}) construction, to the relation  $\mathcal R_G|_F$. The generators of the relation $\mathcal R_G|_F$  are a set of transformations of the quotient space $F\cong \cX/ \Gamma$, in a one to one correspondence with the   cosets of $\Gamma$ in $G$.  We prove that the  composition formula for these transformations is  an averaged version, with coefficients in
$L^\infty(F,\nu)$,   of the Hecke algebra product formula (\cite{BC}).

In the situation $G = \PGL_2(\Z[\frac1p])$, $\Gamma=\PSL_2(\Z)$, $p\geq 3$ prime number,  the relation $\mathcal R_G|_F$ is the equivalence relation associated to a free, measure-preserving  action of a free group on $(p+1)/2$ generators on $F$ (\cite{Ad},\cite {Hj}).  We use the coset representations of the transformations generating $\mathcal R_G|_F$ to find a canonical treeing (\cite{Ga}).

\end{abstract}

\maketitle

\section{Introduction}
Let $G$ be a countable discrete group,
acting ergodicaly by measure-preserving transformations on
an infinite measure space $(\cX,\nu)$ with $\sigma$-finite
measure $\nu$.  We assume that $\Gamma \subseteq G$ is a
discrete  subgroup, that admits a measurable fundamental domain $F$ in $\cX$, of finite measure $\nu(F)=1$, that is
$$\bigcup_{\gamma\in \Gamma} \gamma F=\cX \quad {\rm and} \quad
\nu(\gamma_1 F\cap\gamma_2 F)=0, \quad \gamma_1,\gamma_2\in\Gamma, \gamma_1 \neq \gamma_2.$$

In this paper we study, from an operator algebra point of view,   the probability,   measure-preserving, countable,  measurable, ergodic,  equivalence relation (see e.g. \cite{Ga}) on the probability measure space $(F,\nu|_F)$, defined by the  requirement that two points $x,y\in F$ are equivalent if and only if there exists $g\in G$ such that $gx=y$. We denote this equivalence relation by $\mathcal  R_G|_F$, as this  is the restriction to $F$ of the measure-preserving, countable, measurable, ergodic equivalence relation on $\cX$ determined by the orbits of the action of $G$ on $\cX$ (\cite{Ga}). Such an equivalence relation was also considered in (\cite{Ad}, Example 1.6.3).

Feldman and Moore  constructed (\cite{FM}) a canonical type II von Neumann algebra associated with the relation $\mathcal  R_G|_F$, which we will denote by
$\mathcal M(\mathcal  R_G|_F)$. As explained in the next section, the canonical generators for $\mathcal M(\mathcal  R_G|_F)$ are the elements of the algebra
$L^\infty(F,\nu)$ and a family  of partial measure-preserving isomorphisms (viewed as partial isometries)
$$
\Phi_{\cR_G| F}=\{\phi:A_\phi\rightarrow B_\phi\},$$
where $A_\phi$, $B_\phi$ are measurable subsets of $F$, of equal measure, and the functions $\phi$ are measure-preserving partial isomorphisms (that is isomorphisms on their domain and codomain). We assume that the family $
\Phi_{\cR_G| F}$ is closed under the composition operation (on the maximal possible domain of the composition) and under the inverse operation.

 Since the relation $\mathcal  R_G|_F$ is constructed by restriction to $F$ of the  relation $\mathcal  R_G$, which   has a canonical system of generators - the transformations on $\cX$, given by the elements in the group $G$, there is a canonical choice for the set $\Phi_{\cR_G| F}$
(see Definition \ref{pi} and formula (\ref{choice}) in Lemma \ref{generators}).

 Similarly to the construction of the crossed product algebra, one constructs (\cite{FM}, see also next section) a $\ast$-algebra  $\mathcal M_0(\mathcal  R_G|_F)$ and a positive trace $\tau$ on it. The algebra $\mathcal M_0(\mathcal  R_G|_F)$ is the linear span of $$\Phi_{\cR_G| F}\cdot L^\infty(F,\nu).$$
 The von Neumann algebra $\mathcal M(\mathcal  R_G|_F)$ is constructed from the data consisting of the $\ast$-algebra $\mathcal M_0(\mathcal  R_G|_F)$ and the associated  trace $\tau$, using the standard GNS construction associated to $\tau$.
 By construction, the algebra $\mathcal M_0(\mathcal  R_G|_F)$ has a canonical $\ast$-algebra representation  into $B(L^2(F,\nu))$, which will be described below. We denote this representation by $\piK$, since it is  the restriction of the  Koopman unitary representation (see e.g. \cite{Ke}).

  We denote  by $\mathbb C(G\rtimes L^\infty (\cX, \nu))$  the algebraic crossed product of the algebra $L^\infty (\cX, \nu)$ by the action of the group $G$. By $\mathcal L(G\rtimes L^\infty (\cX, \nu))$ we denote the von Neumann algebra crossed product, obtained by the GNS construction   associated to  the semifinite trace Tr$_\nu$ induced by the $G$-invariant measure $\nu$ on the cross product algebra $\mathbb C(G\rtimes L^\infty (\cX, \nu))$ (see the construction in \cite{FM}).  By $\chi_F\in L^\infty (\cX, \nu) $ we denote the characteristic function of the set $F$. Then
 %%%%%%%%%%
 \begin{equation}\label{corner}
   \mathcal M_0(\mathcal  R_G|_F)=(\chi_F)\big[ \mathbb C(G\rtimes L^\infty (\cX, \nu))\big] (\chi_F).
   \end{equation}
   \noindent and
  \begin{equation}\label{cornervn}
   \mathcal M(\mathcal  R_G|_F)=(\chi_F)\big[ \mathcal L(G\rtimes L^\infty (\cX, \nu))\big] (\chi_F).
   \end{equation}

Here $\chi_F$ is identified with the unit element of the algebra $\mathcal M_0(\mathcal  R_G|_F)$. The canonical trace Tr$_\nu$ (see \cite{FM}), on the crossed product von Neumann  algebra $\mathcal L(G\rtimes L^\infty (\cX, \nu))$, induced from the measure $\nu$ on $\cX$, restricts in the above equality to the trace $\tau$ on $\mathcal M(\mathcal  R_G|_F)$.  The
Koopman $\ast$-algebra representation of
$\mathbb C(G\rtimes L^\infty (\cX, \nu))$ into $B(L^2 (\cX, \nu))$ (see e.g. \cite{Ke}), when restricted to the corner algebra described  in formula
(\ref{corner}), coincides with the representation $\piK$ referred to above.

We use  now the assumption that the group $\Gamma$ is an almost normal subgroup in $G$ (for the definition see formula (\ref{gsigma}) in the  next section). Associated to this data, one constructs the Hecke  $\ast$-algebra $\mathcal H_0(\Gamma,G)$ of  the double cosets of $\Gamma$ in $G$ (for definitions see e.g. \cite{BC} or also formula (\ref{prodhecke}) in  Section \ref{product}). We also have a $\ast$-algebra representation $T=T_{\piK}$ of the algebra $\mathcal H_0(\Gamma,G)$ into $B(L^2 (\cX, \nu))$. The images through $T$ of double cosets in the Hecke algebra are  the Hecke operators. As explained in the next section, the representation $T$ is also associated to the Koopman unitary representation
of $G$ into $B(L^2 (\cX, \nu))$, defined in the next section.

The algebra $\mathcal M_0(\mathcal  R_G|_F)$ considered in this paper is similar to the construction in \cite{RP}, \cite{LLN}. The only difference is that, by making the additional assumption on the existence of a fundamental domain $F$ for the action of the group $\Gamma$ on $\cX$, we have concrete realizations of the quotient space $\cX/\Gamma$ and of this algebra (e.g. as a corner algebra as in formula (\ref{cornervn})). As a bi-product, we obtain an abstract description of the algebra $\mathcal M_0(\mathcal  R_G|_F)$ by generators and relations (see Theorems \ref{formula} and \ref {T}), with generators in
the family $\Phi_{\cR_G| F}$ and the algebra  $L^\infty(F,\nu)$.

The Hecke algebra $\mathcal H_0(\Gamma,G)$ is a subalgebra of $\mathcal M_0(\mathcal  R_G|_F)$. The above Hecke algebra representation $T$
into $B(L^2 (\cX, \nu))$ is obtained by restricting $\piK$ from $\mathcal M_0(\mathcal  R_G|_F)$ to the Hecke algebra $\mathcal H_0(\Gamma,G)$.

When $G = \PGL_2(\Z[\frac1p])$, $\Gamma=\PSL_2(\Z)$, $p$ odd prime number,  the relation $\mathcal R_G| F$ is the equivalence relation associated to a free measure-preserving  action on $F$ of the free group $\mathcal F=F_{(p+1)/2}$ on $(p+1)/2$ generators.
In this situation we construct an explicit treeing (in the sense of \cite{Ga}) for $\mathcal R_G|_F$.
We prove that the  generators and their inverses, which give the treeing, are
 bijective   transformations in the canonical set
 $
\Phi_{\cR_G| F}$ introduced above (constructed in Lemma \ref{generators}).

\section{Main definitions and outline}\label{outline}

In the framework introduced in the previous section,
the subgroup $\Gamma$ is  assumed to be almost normal in $G$. This assumption means that all subgroups
 \begin{equation}\label{gsigma}
\Gamma_\sigma=\sigma\Gamma \sigma^{-1}\cap \Gamma,\quad  \sigma \in G,
\end{equation}
have finite index $[\Gamma:\Gamma_\sigma]$. This property indicates that the subgroup $\Gamma$, which is  not necessary a normal subgroup of $G$,  has, modulo finite (default) indices, properties similar to normal subgroups.
Throughout the paper $\Gamma\sigma$ will denote the left cosets of $\Gamma$ in $G$, while $\Gamma_\sigma$ will denote the subgroup  introduced in (\ref{gsigma}). That formula indicates that the subgroup
$\Gamma_\sigma$  depends only on the left coset $\Gamma \sigma^{-1}$ (or, respectively,  depends only on the right coset $\sigma\Gamma$) for $\sigma\in G$.

Throughout the paper it is assumed that
\begin{equation}\label{sigmaequal}
[\Gamma:\Gamma_\sigma]=[\Gamma: \Gamma_{\sigma^{-1}}], \quad \sigma \in G.
\end{equation}

This will be used  in the construction of the inverses of the partial isomorphisms from the set  $
\Phi_{\cR_G| F}$ (see formula (\ref{choice}) in the Lemma \ref{generators}).
The typical example for a pair of discrete groups as above is (see e.g. \cite{Kr})
 $$\Gamma=\PSL(2, \Z)\subseteq G=PGL(2,\mathbb Q).$$

\noindent The above properties of the pair $\Gamma\subseteq G$ are used to construct a $\ast$-algebra structure on the algebra of double cosets. This algebra is the Hecke algebra $\mathcal H_0(\Gamma, G)$ of double cosets
 of $\Gamma$ in $G$. We refer to \cite{BC}, \cite{Kr}, \cite{RH}, \cite{Tz} for the details of this construction (see also formula (\ref{prodhecke}) in  Section \ref{product}); however we will briefly recall  below the construction of the Hecke operators associated  with a unitary representation of the group $G$.
 The Hecke operators are  then defining the representation of the Hecke algebra (see e.g. \cite{RH}).

   The algebra $\mathcal H_0(\Gamma, G)$ has a canonical  action on the linear space $\mathbb C(\Gamma\backslash G)$ that is freely generated by the coset space $\Gamma\backslash G$. This gives a $\ast$-embedding of the algebra $\mathcal H_0(\Gamma, G)$  into $B(\ell  ^2(\Gamma\backslash G))$. The uniform closure of the image of $\mathcal H_0(\Gamma, G)$ in this embedding is called the reduced Hecke $C^\ast$-algebra of the pair $\Gamma\subseteq G$, and will be denoted by
   $\mathcal H (\Gamma, G)$. Denote by $[\Gamma]$ the vector in $\ell  ^2(\Gamma\backslash G)$ corresponding to the identity coset.  The condition in formula (\ref{sigmaequal}) implies that the state
   $\langle \,\cdot \, [\Gamma], [\Gamma]\rangle$ on $B(\ell  ^2(\Gamma\backslash G))$ restricts to a trace on $\mathcal H (\Gamma, G)$
   (see \cite{BC}).

   Given a representation $\pi$ of the group $G$ into the linear isomorphisms group of a vector space $\mathcal V$, one constructs (see e.g. \cite{RH}, \cite{BC}) an algebra representation $T=T_\pi$ of the Hecke algebra $\mathcal H_0(\Gamma, G)$ into the endomorphism space of the vector space $\mathcal V^{\Gamma}$ consisting of vectors fixed by the subgroup $\Gamma$. Clearly the representation $T$ is determined by its values
   $T([\Gamma\sigma\Gamma])$
   on the double cosets of $\Gamma$ in $G$, $\sigma\in G$, which are a linear generating set
   for the Hecke  algebra $\mathcal H_0(\Gamma, G)$.
The construction  of the linear operators $T([\Gamma\sigma\Gamma])$ is described below. These are the Hecke operators (\cite{He}, \cite{Pe}, see also \cite{Kr}) associated with the representation $\pi$.

By hypothesis, $[\Gamma:\Gamma_\sigma]<\infty $ for every $\sigma \in G$. Let
\begin{equation}\label{cosetdec}
\Gamma=\bigcup_{i=1}^{[\Gamma:\Gamma_\sigma]} s_i [\Gamma_\sigma],\quad s_i\in \Gamma,
\end{equation}
 be the disjoint decomposition of $\Gamma$ as an union of  left cosets of $\Gamma_\sigma$ in $\Gamma$. Since obviously
 $$\sigma [\Gamma_{\sigma^{-1}}]=[\Gamma_\sigma] \sigma,$$
    it follows that every double coset $\Gamma\sigma\Gamma$ as above is a finite disjoint reunion of right cosets of $\Gamma$ in $G$,
    \begin{equation}\label{doublecosets}
\Gamma\sigma\Gamma =\bigcup_{i=1}^{[\Gamma:\Gamma_\sigma]}  s_i\sigma \Gamma ,
\end{equation}
where the coset representatives $s_i\in \Gamma$ are chosen as  in formula (\ref {cosetdec}). Note that a similar, finite decomposition holds true for   left  cosets of $\Gamma$ in $G$.

Then
%for $v$ in $\mathcal V^{\Gamma}$
 the Hecke operators $T=T_\pi$
are defined by the formula
\begin{equation}\label{genhecke}
T([\Gamma\sigma\Gamma])v=\sum_{i=1}^{[\Gamma:\Gamma_\sigma]}
\pi(s_i\sigma)v, \quad  \sigma\in G, v \in \mathcal V^{\Gamma}, [\Gamma\sigma\Gamma]\in \mathcal H_0(\Gamma, G).
\end{equation}
The property that $T([\Gamma\sigma\Gamma])v\in \mathcal V^{\Gamma}$ for $v\in \mathcal V^{\Gamma}$ follows because any $\gamma\in \Gamma$ permutes the cosets in $\Gamma/\Gamma_\sigma$, (see e.g. the argument in \cite{Sa}, page 27). This is indeed a representation of the Hecke algebra (see e.g. \cite{Tz}).

If the vector space $\mathcal V$ contains a $\pi (G)$-invariant Hilbert subspace $H$, and if the family of isomorphisms
$\{\pi(g)|_H, g\in G\}$ is a unitary representation of $G$, then in particular examples one constructs a Hilbert space
$H^\Gamma\subseteq \mathcal V^\Gamma$ (see e.g. \cite{Pe}). In this case the Hecke operators $T([\Gamma\sigma\Gamma])$ leave the Hilbert space $H^\Gamma$ invariant, and their restriction to $H^\Gamma$ extends by linearity to a $\ast$-algebra  representation of the $\ast$-algebra  $\mathcal H_0(\Gamma, G)$.

  The Ramanujan-Petersson conjectures ask for conditions on the representation $\pi$ such that the above representation of the Hecke algebra
$\mathcal H_0(\Gamma, G)$ into $B(H^\Gamma)$ extends to a continuous $C^\ast$-algebra representation of the reduced Hecke $C^\ast$-algebra
$\mathcal H(\Gamma, G)$ introduced above (see \cite{Ra} for this particular formulation of the statement of the Ramanujan-Petersson conjectures). This conjecture has been solved affirmatively, in the context of automorphic forms, by P. Deligne (\cite{De}).

One situation where the Ramanujan-Petersson conjecture is still open is when $\pi$ is a Koopman unitary representation associated with a
measure-preserving transformation of the group $G$ on an infinite measure space $(\cX,\nu)$. We assume, as in the introduction, that the  subgroup $\Gamma$ admits a fundamental domain $F$ in $\cX$ with $\nu(F)=1$. In this case $H=L^2(\cX,\nu)$ and the Koopman unitary  representation $\pi=\piK$ is given by
(see e.g. \cite{Ke})
\begin{equation}\label{koop1}
\piK(g)f(x)=f(g^{-1}x), \quad  x\in \cX,  f\in L^2(\cX,\nu), g\in G.
\end{equation}
In this situation, the larger space $\mathcal V$ is the space of measurable functions on $\cX$. We extend $\pi$ to a representation of $G$ into the linear isomorphism group  of the vector space $\mathcal V$, by using the same formula as above. Clearly the space $\mathcal V^\Gamma$ is identified with the space of measurable functions on $F$. The natural choice for the space $H^\Gamma$ is then $L^2(F,\nu|_F)$. If $\sigma$ is an element in $G$, using the choice of coset representatives as in formula (\ref{doublecosets}),
% for left cosets instead of right cosets of $\Gamma$ in $G$,
%
 the general formula (\ref{genhecke}) gives a Hecke operator
 $$T(\Gamma\sigma\Gamma)= T_{\piK}(\Gamma\sigma\Gamma),$$
 acting  as a bounded linear operator on $L^2(F,\nu|_F)$.
 % We have to switch from right cosets to left cosets because the formula of the Koopman unitary representation at $g\in G$ involves the action of $g^{-1}$.
  The Hecke operator has consequently the formula \begin{equation}\label{hecke}
[T(\Gamma\sigma\Gamma)](f)(x)= \sum_{i=1}^{[\Gamma:\Gamma_\sigma]}
f((s_i\sigma )^{-1}x),\quad  x\in \cX,  f\in L^2(F,\nu), \sigma \in G.
\end{equation}
As explained before, this gives a $\ast$-algebra representation of the Hecke algebra of double cosets $\mathcal H_0(\Gamma, G)$. Since the constant function $1\in L^2(F,\nu)$ is left invariant by all operators (and their adjoints) in formula (\ref{hecke}), it follows that restricting to the orthogonal space $$L_0^2(F,\nu)=L^2(F,\nu)\ominus \mathbb C 1,$$
we obtain a $\ast$-algebra representation of the Hecke algebra $\mathcal H_0(\Gamma, G)$ into $B(L_0^2(F,\nu))$. The Ramanujan-Petersson conjecture then asks for conditions on the action of $G$ on $\cX$ so that this representation extends to a (continuous) $C^\ast$-algebra representation of the reduced Hecke $C^\ast$-algebra $\mathcal H(\Gamma, G)$ into $B(L_0^2(F,\nu))$.

 If the condition in formula (\ref{sigmaequal}) holds true, then the Hecke $C^*$-algebra $\mathcal H(\Gamma, G)$ is commutative  (\cite{Kr}) and the Hecke operators are selfadjoint (\cite{Pe}). The continuity conditions on the representation of the Hecke algebra may be reformulated as  a condition on the location of the spectrum of the Hecke operators  $T(\Gamma\sigma\Gamma)|_{L_0^2(F,\nu)}, \sigma \in G$. In the case $G = \PGL(2, \mathbb Q)$, $\Gamma=\PSL_2(\Z)$, where $\cX$ is the upperhalf complex plane $\mathbb H$ and $G$ acts by M\"obius transformations, these conditions are replaced by   a condition on the point spectrum of the Hecke operators (see e.g. \cite{Sa}).

Recall that the  discrete group $G$ is
acting ergodicaly and a.e. freely, by measure-preserving transformations on
an infinite measure space $(\mathcal X,\nu)$ with $\sigma$-finite
measure $\nu$. We assumed that the restriction of this action to $\Gamma$ admits a fundamental domain $F$ of measure 1.
 We consider the countable  measure-preserving equivalence relation
$\cR_G$ on $\cX$ induced by the orbits of $G$ (see e.g. \cite{Ga}), and its restriction $\cR_G|_F$
to $F$. We will  analyze the  Hecke operators
associated to this data, introduced in formula (\ref{hecke}),  from the point of view of the countable measure-preserving equivalence relation $\cR_G|_F$.

 The standard construction in \cite{FM} associates a type II von Neumann algebra to any measure-preserving, countable, ergodic equivalence relation on a probability space. The construction of the algebra $\mathcal M(\cR_G|_F)$ associated to $\cR_G|_F$
 is outlined below. Generally, because of the selection principle   (\cite{Ku}, see also \cite{FM}), the equivalence relation $\cR_G|_F$ is implemented  by a family  of measure-preserving isomorphisms
\begin{equation}\label{partial}
\Phi_{\cR_G|_F}=\{\phi:A_\phi\rightarrow B_\phi\},
\end{equation}
where $A_\phi$, $B_\phi$ are measurable subsets of $F$, of equal measure, and the maps $\phi$ are measure-preserving isomorphisms, defined on $A_\phi$ and taking values onto $B_\phi$. We refer to such maps as partial isomorphisms.  We assume that $\Phi_{\cR_G|_F}$ is closed with respect to the inverse operation, and that the composition of any two elements  in $\Phi_{\cR_G|_F}$ on the maximal admissible domain is the restriction of another element in $\Phi_{\cR_G|_F}$.   In the particular case of the relation $\cR_G|_F$,  the  set $\Phi_{\cR_G|_F}$ is explicitly described
(see Definition \ref {pi} and Lemma \ref{generators}), being shown to consist of partial isomorphisms, obtained by restricting elements in the group $G$ to suitable measurable subsets of $F$. The Hecke operators are calculated as sums of elements in the set $\Phi_{\cR_G|_F}$.

One considers the partial crossed product algebra $\mathcal M_0(\cR_G|_ F)$, generated as linear space by elements of the
form
\begin{equation}\label{gens}
\phi f, \phi\in \Phi_{\cR_G| F},\quad   f \in L^{\infty}(X,\nu).
\end{equation}

The characteristic function of a measurable subset $A$ is denoted by $\chi_A$.
The product formula on the algebra $\mathcal M_0(\cR_G| F)$ is
determined, for all $\phi\in\Phi_{\cR_G| F}$, $f\in L^{\infty}(X,\nu)$, by  (see \cite{FM})
\begin{equation}\label{crossed}
\phi f= \phi f\chi_{A_\phi}= (f\circ{\phi}^{-1}) \chi_{B_\phi} \phi=(f\circ{\phi}^{-1}).
\end{equation}
In \eqref{crossed} it is implicitly assumed that in the algebra $\mathcal M_0(\cR_G| F)$ the element $\phi$ is identified with  a partial isometry, acting on the Hilbert space $L^2(F,\nu)$, with initial space the projection $\chi_{A_\phi}\in L^{\infty}(X,\nu)$ and range
space equal to the projection $\chi_{B_\phi}\in L^{\infty}(X,\nu)$.
The above formula also determines   the involution on the algebra $\mathcal M_0(\cR_G| F)$, by requiring that
$$(\phi f)^\ast= (\phi)^{-1} \overline {f}, \quad \phi\in \Phi_{\cR_G| F}, f \in L^{\infty}(X,\nu).$$
Here $\overline {f}$ denotes the complex conjugate function.
  Denote by ${\rm Fix}( \phi)$ the subset of fixed points of $\phi$. Since the partial isomorphisms in the family $\Phi_{\cR_G| F}$  are measure-preserving, it follows that the formula
\begin{equation}\label{trace}
\tau_0(\phi f)= \int_{{\rm Fix} (\phi)} f{\rm d}\nu
\end{equation}
defines a trace $\tau_0$ on the algebra $\mathcal M_0(\cR_G| F)$. The arguments in
\cite{FM} prove that, by considering the GNS construction associated to $\tau_0$, one obtains
a type II von Neumann algebra $\mathcal M(\cR_G|_F)$ endowed with a finite  faithful trace $\tau$. If the measure-preserving relation $\cR_G|_F$ is ergodic (see e.g. \cite{Ga}), then
$\mathcal M(\cR_G|_F)$ is a factor (\cite{FM}).  We let $\mathcal A_{\rm red}(\cR_G|_F)$ be the norm closure
of the $\ast$-algebra $\mathcal M_0(\cR_G| F)$ in $\mathcal M(\cR_G| F)$. This will be referred to as the reduced
$C^\ast$-algebra associated to the relation $\cR_G|_F$.

The involutive algebra $\mathcal M_0(\cR_G| F)$ has a canonical $\ast$-algebra representation $\piK$, the analogue of the Koopman unitary representation
into the bounded linear operators on $L^2(F,\nu)$.
This is obtained as follows: Let $M_f\in B (L^2(F,\nu))$ denote the multiplication operator by $f\in \L^\infty(F,\nu)$. For a partial isomorphism $\phi\in \Phi_{\cR_G| F}$ as in formula, denote by $U_\phi$ the partial isometry on $\L^\infty(F,\nu)$ with initial, respectively final, projection $M_{\chi{A_\phi }}$, respectively $M_{\chi{B_\phi }}$, defined by mapping $f\in L^2(A_\phi,\nu)$ into
$L^2(A_\phi,\nu)$.
The precise formula for the partial isometry $U_\phi$ is
\begin{equation}\label{uphi}
U_\phi(f)(x)=f(\phi^{-1}x), \quad x \in A_\phi, f\in L^2(B_\phi,\nu).
\end{equation}

 The representation $\piK$  is defined on generators   as follows:
\begin{equation}\label{koop}
\piK(f)= M_f;\quad  \piK(\phi)=U_\phi, \quad f \in \L^\infty(F,\nu), \phi\in \Phi_{\cR_G|_F}.
\end{equation}
The above formula extends by linearity to a $\ast$-representation $\piK$  of the $\ast$-algebra $\mathcal M_0(\cR_G|_F)$. 

 The representation introduced in formula (\ref {koop1}) has a canonical extension, denoted by $\tilde{\piK}$, to the algebra $\mathbb C(G\rtimes L^\infty (\cX, \nu))$, with values in $B(L^2(\cX,\nu))$. As shown in formula (\ref{corner}), the  algebra $\mathcal M_0(\cR_G|_F)$ is a reduced  algebra  of the algebra $\mathbb C(G\rtimes L^\infty (\cX, \nu))$.  The representation $\piK$ introduced  in the  formula (\ref{koop})  is the restriction of the representation $\tilde{\piK}$ to the algebra $\mathcal M_0(\cR_G|_F)$. 

To keep notation simple, when considering a partial isomorphism $\phi$ as a partial isometry $U_\phi$ in the algebra $M_0(\cR_G|_F)$, we will denote $U_\phi$ by $\phi$. Then, the  composition of two  partial isomorphisms
$\phi,\psi$ on the maximal admissible domain of the composition operation  corresponds to the product
$$ U_\phi U_\psi\in M_0(\cR_G|_F) .$$

More generally, let
$$\dot{\wideparen{H}}:F\rightarrow F$$
be a measurable transformation (not necessary bijective).  Assume that there exists two measurable partitions of the set $F$,  $(A_i)_{i=1}^n$,  $(B_i)_{i=1}^n$, such that $\dot{\wideparen{H}}$ maps $A_i$ onto $B_i$, and
$$\dot{\wideparen{H}}|_{A_i} :A_i\rightarrow B_i,\quad  i=1,2,\ldots,n,$$
is a partial, measure-preserving isomorphism, as above. In particular the sets $A_i$, $B_i$ have equal measure for all $i$.  We will refer to the sets
$A_1,\ldots,A_n$ as domains of injectivity (or bijectivity) for the transformation $\dot{\wideparen{H}}$, writing this formally as
\begin{equation}\label{domains}
\dot{\wideparen{H}}=\sum_{i=1}^n [\dot{\wideparen{H}}]\chi_{A_i}=\sum_{i=1}^n \chi_{B_i}[\dot{\wideparen{H}}].
\end{equation}
To each element as above we associate in the algebra $\mathcal M_0(\cR_G|_F)$ the element
\begin{equation}\label{uh}
\sum_{i=1}^n U_{\dot{\wideparen{H}}}\chi_{A_i}.
\end{equation}
The notational convention introduced above implies that formula  (\ref{domains}) remains
valid, if we consider the transformation $\dot{\wideparen{H}}$ as an element of the algebra
$M_0(\cR_G|_F)$. Clearly, by linearity, the composition $$
\dot{\wideparen{H_1}}\circ \dot{\wideparen{H_2}},$$
of two transformations $\dot{\wideparen{H_1}}$, $\dot{\wideparen{H_2}}$ as above corresponds, in the setting introduced here, to the product of the two elements
associated to the transformations, $\dot{\wideparen{H_1}}$, $\dot{\wideparen{H_2}}$ in the algebra $M_0(\cR_G|_F)$.

The $C^\ast$-algebra generated by the image through the representation $\piK$ of $\mathcal M_0(\cR_G|_ F)$,
introduced in formula (\ref{koop}), will be denoted by $\mathcal A_{\operatorname{Koop}}(\cR_G|_F)$, that is
\begin{equation}\label{akoop}
\mathcal A_{\operatorname{Koop}}(\cR_G|_F)=\overline{ \piK (\mathcal M_0(\cR_G|_F))}^{||\cdot||}
\subseteq B(L^2(F)).
\end{equation}

The representation $T=T_{\piK}$ introduced in formula (\ref{hecke}) gives a $\ast$- algebra embedding of the Hecke algebra into
$\mathcal A_{\operatorname{Koop}}(\cR_G| F)$. Hence, the determination of the
$C^\ast$-norm on the algebra  $\mathcal A_{\operatorname{Koop}}(\cR_G|_F)$ is essential for the determination of the continuity properties of the Hecke operators.

As specified above, we will also use the notation $\phi$ for the unitary $U_\phi$ when $\phi \in \Phi_{\cR_G|_F}$. The almost normality of $\Gamma$ in $G$ is used in Theorem \ref{formula} to prove that the composition rule of the elements in the set  $\Phi_{\cR_G| F}$ has a formula analogue to the Hecke algebra multiplication,   with $L^\infty(F,\nu)$-coefficients.
    In Theorem \ref{T} we establish a precise presentation by generators and relations for
the  algebra $\mathcal M(\cR_G|_F)$.
This is an  explicit description for the action of the generators of
$\cR_G|_F$  which is context free (does not depend on $F$)
in the sense of symbolic dynamics.

In the situation $G = \PGL_2(\Z[\frac1p])$, $\Gamma=\PSL_2(\Z)$, $p\geq 3$ prime number,
 the formula in Theorem \ref{formula} for the  composition
of the generators of the $*$-algebra associated to the equivalence relation
$\cR_G|_F$,
proves that
 $\cR_G| F$ is treeable and
of cost $\frac{p+1}{2}$.
% By the results
%of Hjorth, the equivalence relation is therefore implemented by the action of
%the free group with $\frac{p+1}{2}$ generators on $F$.
In Theorem \ref{prime}we prove that the generators introduced in Definition \ref{pi} and Lemma \ref{generators} provide a canonical  treeing  (\cite{Ga})  for the treeable relation $\cR_G|_F$.
%Moreover, the radial algebra of the free group (the algebra generated by
%convolutors in the words on $F_{\frac{p+1}{2}}$ of equal length
%have equal weight) will coincide with the Hecke algebra corresponding to $G$,
%$\Gamma$ and to the action on~$\cX$.

In particular, it follows as a corollary that, in analogy with the notion of measured equivalence for
groups (\cite{Ga}), the group $\PGL_2(\Z [\frac1p])$ is infinitesimally
orbit equivalent to $F_{(p+1)/2}$ (see Corollary \ref{ior} for the definition of infinitesimally
orbit equivalence).

\section{Generators of the  relation $\cR_G|_F$}\label{product}

In this section we construct a canonical   family of generators, as in formula (\ref{partial}), for the relation $\cR_G|_F$.  We analyze first the composition formula for such generators, and prove that this formula is a
$L^\infty(F,\nu)$-coefficients variant of the Hecke algebra multiplication formula for double cosets. We determine the formula for the inverses
of these transformations on bijectivity domains.
% and determine a canonical  set of partial isomorphisms $\Phi_{\cR_G|_F}$ generating the equivalence relation.

\begin{defn}\label{pi}
 Let $G$ be a discrete group acting by measure-preserving
transformations, almost everywhere freely on  $\cX$. Assume that $\Gamma$ is an almost normal subgroup, having a  fundamental domain $F\subseteq \cX$ of finite measure.

For $g$ in $G$, we define a (non-injective) transformation (function)
$$\dot{\wideparen{\Gamma g}}:F\rightarrow F,$$ as follows:

Let $x$ be an element  in $F$. Since $F$ is a fundamental domain, there exists
a unique $\gamma_1\in \Gamma$ and $x_1$ in $F$ such that $gx=\gamma_1x_1$. Then, we define:
\begin{equation}\label{gammag}
\dot{\wideparen{\Gamma g}} (x):=x_1=\gamma_1^{-1}gx.
\end{equation}
 Clearly, the function $\dot{\wideparen{\Gamma g}}$
depends only on the right $\Gamma$-coset $\Gamma g$  of $g\in G$.

Consequently $\cR_G|_F$ is generated by the transformations $\dot{\wideparen{\Gamma g}}$, with $g$ running
through a system of representatives for right cosets of $\Gamma$.
Hence, for $x,y\in F$ we have that  $x\sim y$ with respect to the equivalence relation $\cR_G|_F$  if and only if there exists $g\in G$ such that
$\dot{\wideparen{\Gamma g}} x=y$.

In Lemma \ref{bijective} below we will prove that the transformations $\dot{\wideparen{\Gamma g}}$ are of the type considered in formula (\ref{domains}).
\end{defn}

\begin{lemma}
For every $g\in G$, the transformation $\dot{\wideparen{\Gamma g}}$ is not injective, but the cardinality
of the preimage of each point in the image is equal to $[\Gamma:\Gamma_g]$.
In addition, if $\Gamma gs_i$ are the left $\Gamma$-cosets contained in   $\Gamma g \Gamma$,
then every point $x$ in $F$ will show up exactly $[\Gamma:\Gamma_g]$-times in the
reunion of the images of the maps $\dot{\wideparen{\Gamma g}}s_i$.

The same
is true for preimages, with $[\Gamma:\Gamma_{g^{-1}}]$ instead of $[\Gamma:\Gamma_g]$.

\end{lemma}
%\begin{proof}
%To prove this, we note that an alternative definition, for the transformations $\dot{\wideparen{\Gamma g}}$, is as follows.
%If $f$ is a positive valued, measurable function on $\cX$, we define
%$$\Gamma f=: \sum_{\gamma\in Gamma} \pi_{Koop}(\gamma) f.$$
%Then $\Gamma f$ is a $\Gamma$-invariant positive function (with possible infinite values).
%\end{proof}

\begin{proof} We have to  count images and preimages.
But this follows from the fact that if $s_i$ are the coset representatives for $\Gamma_{g^{-1}}$ in $\Gamma$, then the domain $F_0=\bigcup_i g s_i F$ is a fundamental
domain for $\Gamma_{g^{-1}}$ and it covers $[\Gamma:\Gamma_g]$ times the set $F$.
%Here $s_i$ are the right  coset representatives for the subgroup $\Gamma g$ of $\Gamma$.
\end{proof}

The transformations $\dot{\wideparen{\Gamma g}}$, $g\in G$ have a natural composition rule,
 similar to the multiplication rules from the Hecke
algebra with $L^\infty(F,\nu)$ coefficients.

\begin{thm}\label{formula}
With the previous hypothesis, let $g_1,g_2$ be arbitrary elements in $G$.
Let $r_j$, $j=1,2,\ldots,[\Gamma:\Gamma_{g_1^{-1}}]$
be a family of right coset representatives for $\Gamma_{g_1^{-1}}$ in $\Gamma$.
Let $A_{g_1,g_2}^{r_j}$ be the subset of $F$ defined by the formula:
$$
\{f\in F\mid r_jg_2f\in \Gamma_{g_1^{-1}}F\}=(r_jg_2)^{-1}\Gamma_{g_1^{-1}}F\cap F.
$$
Then
\item (i).
 The composition of the  transformations  $\dot{\wideparen{\Gamma g_1}}$ and $\dot{\wideparen{\Gamma g_2}}$ of $F$ is given by \begin{equation}\label{mixedhecke}
\dot{\wideparen{\Gamma g_1}}\circ \dot{\wideparen{\Gamma g_2}}=
\sum_j \dot{\wideparen{\Gamma g_1r_j g_2}} \ \chi_{A^{r_j}_{g_1,g_2}}.
\end{equation}

%Here $\dot{\wideparen{\Gamma g_1}}\circ \dot{\wideparen{\Gamma g_2}}$ corresponds to the composition of the two transformations $\dot{\wideparen{\Gamma g_1}}\; $ $\dot{\wideparen{\Gamma g_2}}.$
%

\item (ii) The sets ${A_{g_1,g_2}^{r_j}}$, $j=1,2,\ldots, [\Gamma:\Gamma_{g^{-1}}]$
give a partition of the set $F$.
\end{thm}

Formula (\ref{mixedhecke}) indicates that the term  on the  left side of the equation,  when restricted
to any of the sets ${A_{g_1,g_2}^{r_j}}$, for a fixed index $j$, is equal to
the transformation $\dot{\wideparen{\Gamma g_1r_j g_2}}$ restricted to the same set.

Before giving the proof of the theorem, we note the similarity between formula (\ref{mixedhecke}) and the representation (\cite{BC}, \cite{Tz}) of the  Hecke algebra into $B(\ell^2(\Gamma\backslash G))$.
Indeed, using the notation from the statement of the above theorem, we have that the following formula defines the action of the double coset $[\Gamma g_1\Gamma]\in \mathcal H_0(\Gamma, G)$ on the right coset $[\Gamma g_2]\in \ell^2(\Gamma\backslash G)$:
\begin{equation}\label{prodhecke}
[\Gamma g_1\Gamma][\Gamma g_2]=\sum_j [\Gamma g_1r_j g_2].
\end {equation}
This formula, using the fact that any double coset is a finite union of left (respectively right) cosets, determines the product in the Hecke algebra
(see e.g. \cite{Kr}).

\vskip6pt
\begin{proof}[Proof of Theorem \ref{formula}] Since $G$ acts freely almost everywhere, we may simply work on an
orbit of $G$. So we way assume that $X=G$, and that $F=S\subseteq G$ is a system of coset representatives for
$\Gamma\setminus G$.
Our assumption means that $\Gamma=\bigcup_j \Gamma_{g_1^{-1}} r_j$.  Equivalently we have the disjoint union $\Gamma g_1\Gamma=\bigcup_j \Gamma g_1r_j$.

%(The only point where the initial data would enter
%would be in the measure of the sets in the partitions from the previous proposition.)

Given $s\in S$ and two left cosets $\Gamma g_1,\Gamma g_2$,
we calculate the composition $$[\dot{\wideparen{\Gamma g_1}}]\circ[\dot{\wideparen{\Gamma g_2}}] (s).$$

\noindent We assume that $g_2s=\gamma_2s_2$ for some $\gamma_2\in\Gamma$, $s_2\in S$,
and thus $$\dot{\wideparen{\Gamma g_2}}(s)=s_2.$$

\noindent Then $s_2=\gamma_2^{-1}g_2s$ and hence
$$
g_1([\dot{\wideparen{\Gamma g_2}}](s))=g_1s_2=g_1\gamma_2^{-1}g_2s.
$$
We first identify  the coset of $\Gamma_{g_1^{-1}}$, to which the element $\gamma_2^{-1}$ belongs.
Assume thus that $\gamma_2^{-1}$ belongs to $\Gamma_{g_1^{-1}}r_j$ for some fixed $j$.
Thus we are assuming that  $\gamma_2^{-1}=\theta r_j$ for some $\theta\in \Gamma_{g_1^{-1}}$.
Then $$g_1\gamma_2^{-1} g_2s=(\gamma_1\theta\gamma_1^{-1})g_1r_jg_2s.$$
But $\theta'=\gamma_1\theta\gamma_1^{-1}$ belongs to $\Gamma_g\subseteq\Gamma$, since
$$g\Gamma_{g^{-1}}g=g(\Gamma\cap g^{-1}\Gamma g)g= \Gamma_g.$$
Thus
$$
g_1(\dot{\wideparen{\Gamma g_2}}s)=g_1\gamma_2^{-1}g_2s=\theta'(g_1r_jg_2)s.
$$
On the other hand, there exist $\gamma_1\in\Gamma$, $s_1\in S$ such that $$g_1r_jg_2s=\gamma_1s_1.$$
Thus $$\dot{\wideparen{\Gamma g_1 r_jg_2}}(s)=s_1.$$
 From the above formula we conclude that
$$
g_1([\dot{\wideparen{\Gamma g_2}}](s))=\theta'\gamma_1s_1 ,
$$
and hence
\begin{equation}\label{composition}
\dot{\wideparen{\Gamma g_1}}\big([\dot{\wideparen{\Gamma g_2}}](s)\big)=s_1=\dot{\wideparen{\Gamma g_1r_jg_2}}(s).
\end{equation}
We have to determine the conditions that we have to impose on $s$, so that $\gamma_2$
belongs to $\Gamma_{g_1^{-1}}r_j$.
But the relation defining $s_2$ was
$$
g_2s=\gamma_2 s_2.
$$

Thus, for $\gamma_2^{-1}$ to belong to  $\Gamma_{g_1^{-1}}r_j$, which is equivalent to
$$\gamma_2\in r_j^{-1} \Gamma_{g_1^{-1}},$$ it is necessary and sufficient that
$g_2s$ belongs to $r_j^{-1} \Gamma_{g_1^{-1}}S$.

Thus, for the given choice of the coset representative $r_j$, formula
(\ref{composition}) holds true if  the element $s$ belongs to the set
 $$A_{g_1,g_2}^{r_j}=g_2^{-1}r_j^{-1}\Gamma_{g_1^{-1}}S\cap S.$$
Hence, formula (\ref{mixedhecke})  holds true  on the set $A_{g_1,g_2}^{r_j}$. This completes the proof of (i).

Part (ii) is proved as follows.
Since the cosets $r_j^{-1}\Gamma_{g_1^{-1}}$ are disjoint and $S$ is
a set representatives for $\Gamma\backslash G$, it follows that $\gamma S\cap S=\emptyset$ for all $\gamma\neq e$
and hence $\gamma_1S\cap \gamma_2 S=\emptyset$ if $\gamma_1\neq \gamma_2$. Hence
$$r_j^{-1}\Gamma_{g_1^{-1}}S\cap r_k^{-1}\Gamma_{g_1^{-1}}S =\emptyset \quad \mbox{\rm if $j\neq k$.}$$
Since $$\bigcup_j A_{g_1,g_2}^{r_j}=g\Gamma S\cap S=gGS\cap S=G\cap S=S,$$
 it follows that the sets $A_{g_1,g_2}^{r_j}$, $j=1,2,\ldots,[\Gamma:\Gamma_{g_1^{-1}}]$
give a partition of the set $S$ (and hence of the set $F$, with the initial notation).
\end{proof}

In the context of the previous theorem, we note that the decomposition
$$
\dot{\wideparen{\Gamma g_1}}\circ\dot{\wideparen{\Gamma g_2}}=
\sum_j\dot{\wideparen{\Gamma g_1r_jg_2}}\
\chi_{[g_2^{-1}r_j^{-1}[\Gamma_{g_1^{-1}}]F\cap F]}
$$
obviously depends only on the coset class
$\Gamma g_1$ of the group element $g_1$, since $$\Gamma_{g_1^{-1}}=g_1^{-1}\Gamma g_1\cap\Gamma.$$

The formula does not depend either on the choice of the representative $g_2$
in $\Gamma g_2$, since changing $g_2$ into $\gamma'g_2$ would have
the effect of permuting the sum.  This is due to the fact that
$$
[\Gamma_{g_1^{-1}}] r_j\gamma'= [\Gamma_{g_1^{-1}}]r_{\pi_{\gamma'}(j)}
$$
for a  permutation $\pi_{\gamma'}$ of the set $\{1,2,\ldots,[\Gamma:\Gamma_{g_1^{-1}}]\}$, depending only on  $\gamma'$.

Using the results from the previous statement, we  also prove the following formula for inverses of the transformations introduced in Definition \ref{pi}, when restricted to their corresponding injectivity domains.

\begin{lemma}\label{bijective}
Fix $g\in G$. Let $\alpha_i$, $i=1,2,\dots, [\Gamma:\Gamma_g]$ be
a system of right representatives for cosets of $\Gamma_g$ in $\Gamma$.  Then
\item (i)
For each $i$, the image through $\dot{\wideparen{\Gamma g}}$ of the set
$$A_{\alpha_i,\Gamma g}=:g^{-1}[\Gamma_g]\alpha_i F\cap F=\{s\in F\mid gs\in\Gamma_g\alpha_i F\},$$ is the set
$$B_{\alpha_i,\Gamma g}=:\alpha_i^{-1}[\Gamma_g] gF\cap F.$$

\noindent The transformation $\dot{\wideparen{\Gamma g}}|_{A_{\alpha_i,\Gamma g}}$, with values in the set $B_{\alpha_i,\Gamma g}$,
is bijective. Its inverse is the restriction of the transformation $\dot{\wideparen{\Gamma g^{-1}\alpha_i}}$  to the set
$B_{\alpha_i,\Gamma g}$.

 \item (ii)  The sets $A_{\alpha_i,\Gamma g}$, $i=1,2,\dots, [\Gamma:\Gamma_g]$ give a partition of $F$.
In general, the sets $B_{\alpha_i,\Gamma g}$, $i=1,2,\dots, [\Gamma:\Gamma_g]$ do not form a partition, as they may have overlaps in $F$.
\end{lemma}

\vskip6pt
\begin{proof}
As in the proof of the previous lemma,  we way assume that $X=G$, and that $F=S$ is a system of coset representatives for
$\Gamma\setminus G$.

The choice of the right coset representatives $\alpha_i$, corresponds to the fact that the
 following disjoint unions hold true:
\begin{equation}\label{alpha}
\Gamma=\bigcup_i\Gamma_g\alpha_i \quad \mbox{and} \quad \Gamma g^{-1}\Gamma=\bigcup_i\Gamma g^{-1}\alpha_i.
 \end{equation}

First we show that the image through $\dot{\wideparen{\Gamma g}}$ of
the set $A_{\alpha_i,\Gamma g}$ is the set $\alpha_i^{-1}\Gamma_g  S\cap S$.
Let
$$ s\in A_{\alpha_i,\Gamma g}=g^{-1}\Gamma_g\alpha_i^{-1}  S\cap S=
\{s\in S\mid gs\in \Gamma_g\alpha_i S\}.
$$
Thus $$gs=\theta\alpha_i s_1$$
for some $s_1\in s$, $\theta\in \Gamma g$.
Then $$s_1=\alpha_i=\theta g,$$ and this   belongs to $\alpha_i^{-1}\Gamma_g g S\cap S$.

To verify the inverse formula we have to calculate the composition
$$\dot{\wideparen{\Gamma g^{-1}\alpha_i}}\circ\dot{\wideparen{\Gamma g}}.$$
Using  the previous theorem,  we  let $r_j$, $j=1,2,\dots, [\Gamma:\Gamma_{(g^{-1}\alpha_i)^{-1}}]$ be a system of right representatives
for  the subgroup $\Gamma_{(g^{-1}\alpha_i)^{-1}}$ in $\Gamma$.  This corresponds to the fact that
$$
\Gamma=\bigcup_j [\Gamma_{(g^{-1}\alpha_i)^{-1}}]r_j.
$$
Clearly
 $$\Gamma_{(g^{-1}\alpha_i)^{-1}}=\Gamma_{\alpha_i^{-1} g}=\alpha_i^{-1}\Gamma_g\alpha_i.$$
\noindent Hence we have the disjoint union decomposition
 $$\Gamma=\bigcup_j(\alpha_i^{-1}\Gamma_g\alpha_i)r_j.$$

Because of the initial  assumption from formula (\ref{sigmaequal}), it follows that the two sets of representatives
$(\alpha_i)_i$ and $(r_j)_j$ have the same cardinality.
Then, formula (\ref{mixedhecke}) from the statement of the previous theorem proves  that
\begin{equation}\label{idt}
\dot{\wideparen{\Gamma g^{-1}\alpha_i}}\circ\dot{\wideparen{\Gamma g}}
=\sum_j \dot{\wideparen{\Gamma g^{-1}\alpha_ir_j g}}\
\chi_{[g^{-1}r_j^{-1}[\Gamma_{\alpha_i^{-1}g}]S\cap S]}.
\end{equation}
The identity transformation will appear in the right hand side of the sum above exactly when $\alpha_i r_j$
belongs to $\Gamma_g$.
Consequently, the identity term in the right hand side sum of (\ref{idt})  is located on the set
\begin{equation}\label{inversedomain}
g^{-1}r_j^{-1}[\Gamma_{\alpha_i^{-1}g}]S\cap S=g^{-1}r_j^{-1}[\alpha_i^{-1}\Gamma_g\alpha_i]S\cap S
\end{equation}
iff $\alpha_i r_j$ belongs to $\Gamma_g$. In this case the set in formula
(\ref{inversedomain}) coincides with
$$g(\alpha_i r_j)^{-1}[\Gamma_g]\alpha_i S\cap S=g^{-1}[\Gamma_g]\alpha_i S\cap S.$$
Consequently, the inverse of the transformation $\dot{\wideparen{\Gamma g}}$, when restricted to the domain $g^{-1}[\Gamma_g]\alpha_i S\cap S$,
is the transformation $\dot{\wideparen{\Gamma g^{-1}\alpha_i}}$.
This completes the proof of (i).

%It is easy to see that the inverse formula is consistent:  if
%we apply formula this to $\dot{\wideparen{\Gamma g^{-1}\alpha_i}}$ on $\alpha_i\Gamma_g gS\cap S$
%we get the same result.

Part (ii) is proved as   in the previous lemma, using the fact that $F$ is a fundamental domain for the action of the group $\Gamma$ on $\cX$.
\end{proof}

In the following lemma we list the inverses of all transformations corresponding to the right cosets of $\Gamma$ that belong to a fixed double coset of $\Gamma$. This allows us to list the set of transformations belonging to  an explicit set of transformations $\Phi_{\cR_G|_F}$ as in formula (\ref{partial}), generating the countable equivalence relation $\cR_G|_F$.

For any $g$ in $G$, let $r_j$, $j=1,2,\dots ,[\Gamma:\Gamma_{g^{-1}}]$ be a system of left coset representatives for
the subgroup $\Gamma_{g^{-1}}$ of $\Gamma$. This choice corresponds to
\begin{equation}\label{r}
\Gamma=\bigcup_j \Gamma_{g^{-1}}r_j\quad \mbox{\rm or}\quad
\Gamma g\Gamma=\bigcup_j\Gamma g r_j.
\end{equation}

The inverses of all transformations  $\dot{\wideparen{\Gamma g r_j}}$, restricted to their corresponding bijectivity domains
 $A_{\alpha_i,\Gamma gr_j}$ constructed in the previous lemma, are listed below in formula (\ref{choice}).

   Let the coset representatives
  $(\alpha_i)_i$  be defined by the requirement from formula (\ref{alpha}) in the statement of the previous lemma.
   The initial  assumption expressed in formula (\ref{sigmaequal}) implies that the two sets of coset representatives introduced in formulae (\ref{alpha}) and (\ref{r}) have the same cardinality.
%where
%$r_j$ are a system of left coset representatives for
%$\Gamma_{g^{-1}}$ in $\Gamma$ (thus $\Gamma=\bigcup\Gamma_{g^{-1}}r_j)$
%and then $\Gamma g\Gamma=\bigcup\Gamma gr_j$.

Since
$\Gamma_{g r_j}=\Gamma_g$, the  choice of the coset representatives
$(\alpha_i)_i$ is independent of the choice of the coset representatives $(r_j)_j$ introduced in formula (\ref{r}). The choices of the coset representatives may be summarized by requiring that both disjoint unions
 $\bigcup_i\Gamma_g\alpha_i$ and $\bigcup_j\Gamma_{g^{-1}}r_j$ are equal to $\Gamma$.

\begin{lemma}\label{generators} Let $g\in G $ and let $(\alpha_i)$  and $(r_j)$ denote
the  pair of coset representatives constructed in formulae (\ref{alpha}) and (\ref{r}), with $i,j\in \{1,\dots, [\Gamma:\Gamma_{g^{-1}}]=[\Gamma:\Gamma_{g}]\}$.
Then:
\item (i) For all $i,j$, the inverse of the transformation
$\dot{\wideparen{\Gamma gr_j}}$, restricted to the bijectivity domain
\begin{equation}
\big[(gr_j)^{-1}[\Gamma_{gr_j}]\alpha_i^{-1}\big]F\cap F=\big[(gr_j)^{-1}[\Gamma_g]\alpha_i^{-1}\big]F\cap F,
\end{equation}
is the transformation
$
\dot{\wideparen{\Gamma r_j^{-1}g^{-1}\alpha_i}}=\dot{\wideparen{\Gamma g^{-1}\alpha_i}},
$
restricted to the bijectivity domain
$$
\big[\alpha_i^{-1}[\Gamma_g] gr_j \big]F\cap F=\big[\alpha_i^{-1}g[\Gamma_{g^{-1}}]r_j\big]F\cap F.
$$
\item (ii) The canonical choice for the set $\Phi_{\cR_G|_F}$ as  in formula (\ref{partial}) is given by the following set of  measure-preserving, partial isomorphisms of the set $F$:
\begin{equation}\label{choice}
%\Phi_{\cR_G|_F}=
\big\{ \dot{\wideparen{\Gamma g r_j}}:{\big[r_j^{-1}g^{-1}[\Gamma_g]\alpha_iF\cap F\big]}\rightarrow \big[\alpha_i^{-1}g[\Gamma_{g^{-1}}]r_j\big]F\cap F
\ | \  g\in G \big\},
\end{equation}
where for each double coset   $ \Gamma g \Gamma $, the coset representatives $(\alpha_i,r_j)_{i,j=1,\ldots,[\Gamma :\Gamma_g]}$  are chosen as above.
\end{lemma}
\begin{proof} This follows from the previous lemma.
One applies  the formula for the inverse of  the transformation $\dot{\wideparen{\Gamma g}}$ for the transformation $\dot{\wideparen{\Gamma gr_j}}$
\end{proof}

With the notation from Lemma \ref{generators}, we note that because $F$ is a fundamental domain for the action of the group $\Gamma$, the sets $r_j^{-1}g^{-1}[\Gamma_g]\alpha_iF\cap F$ are disjoint after $i$. Similarly, the sets
$[\alpha_i^{-1}g[\Gamma_{g^{-1}}]r_j]F\cap F$ are disjoint after $j$.

\begin{rem}\label{rem6}
In the context of the previous two  lemmas, fix $g\in G$ and consider the (non-injective) transformation $\dot{\wideparen{\Gamma g}}:F\rightarrow F$.
Using the notation in Lemma \ref{bijective}, choosing $(\alpha_i)$, $i=1,\dots, [\Gamma:\Gamma_g]$ as in formula (\ref{alpha}), we have
\begin{equation}\label{explanation}
\dot{\wideparen{\Gamma g}}=\sum_i \dot{\wideparen{\Gamma g}}
\chi_{[A_{\alpha_i,\Gamma g}]}.
\end{equation}
In the algebra $\mathcal M_0(\cR_G|_F)$ this corresponds, by formula
(\ref{uphi}), to  the element:
\begin{equation}\label{element}
\sum_i U_{\dot{\wideparen{\Gamma g}}|_{A_{\alpha_i,\Gamma g}}}.
\end{equation}
This is possible, as explained in Section \ref{outline} in formulae (\ref{domains}) and (\ref{uh}), because the restriction $\dot{\wideparen{\Gamma g}}|_{A_{\alpha_i,\Gamma g}}$ is injective.

In Section \ref{outline}, when describing the generators for the algebra $\mathcal M_0(\cR_G|_F)$, we introduced the  convention that, as element of the above mentioned algebra, for a measure-preserving partial isomorphism $\phi\in  \Phi_{\cR_G|_F}$, the partial isometry $U_\phi$ in formula (\ref{uphi}) will be simply denoted by $\phi$.
With this notation to the transformation $\dot{\wideparen{\Gamma g}}$ corresponds the  element in the right hand side of formula (\ref{explanation}), when regarded as an element of the $\ast$-algebra $\mathcal M_0(\cR_G|_F)$.

Let $\piK$ be the Koopman representation of $G$ on $L^2(\mathcal X, \nu)$ and let $P_{\chi_F}$ be the projection operator of
multiplication with the characteristic function  $\chi_F$ on $L^2(\mathcal X, \nu)$.
Then, the construction in Definition \ref{pi} proves that the elements in the algebra $\mathcal M_0(\cR_G|_F)$, introduced in the formulae (\ref{explanation}) and (\ref{element}),
also  have the alternative formula
\begin{equation}\label{identif}
\dot{\wideparen{\Gamma g}}=\sum_{\theta \in \Gamma g}P_{\chi_F}\piK (\theta)P_{\chi_F},\quad  g\in G.
\end{equation}
The sum in the right hand side of the above formula is considered in \cite{Ra1}.

Fix a double coset $[\Gamma g\Gamma]$ of $\Gamma$ in $G$ and let
$\mathcal A_{[\Gamma g\Gamma]}$ be the set of partial isomorphisms, enumerated in formula (\ref{choice}), corresponding to a fixed $g\in G$.

Using the results in \cite{Ra1}, we obtain that the Hecke operator associated to the double coset
$$[\Gamma g\Gamma]=\bigcup_i \Gamma\sigma s_i,$$
can be expressed as
\begin{equation}\label{heckeformula}
T([\Gamma g\Gamma])= \sum_{\theta \in \Gamma g\Gamma}P_{\chi_F}\piK (\theta)P_{\chi_F}= \sum_i  \dot{\wideparen{\Gamma gs_i}}
\in \mathcal M_0(\cR_G|_F).
\end{equation}
Thus, using the results and notation from  Lemma \ref{generators}, we obtain the following expression for the Hecke operators, associated to a double coset $[\Gamma g\Gamma]$, $g\in G$:
\begin{equation}\label{heckeformula1}
T([\Gamma g\Gamma])=\sum_{\phi\in \mathcal A_{[\Gamma g\Gamma]}} \phi.
\end{equation}

%Note that the generating  transformations from Definition \ref{pi}, viewed as in formula (\ref{element}) as elements in $\mathcal M_0(\cR_G|_F)$, have also the property that the pairing
%\begin{equation}\label{mult}
%\dot{\wideparen{\Gamma g_1}}, \dot{\wideparen{\Gamma g_2}}\rightarrow
 %[\dot{\wideparen{\Gamma g_2}}]^\ast \cdot \dot{\wideparen{\Gamma g_1}}\in \mathcal %M_0(\cR_G|_F),\quad   g_1,g_2\in G,
% \end{equation}
%is a representation of the obvious pairing
%$$\C( G\backslash \Gamma) \times \C(\Gamma/G)\rightarrow \C(g_1 \Gamma g_2\mid g_1,g_2\in G).$$
%In the above formula the transformations $\dot{\wideparen{\Gamma g_1}}, \dot{\wideparen{\Gamma g_2}}$, $g_1,g_2\in G$ are identified with the corresponding elements in the algebra $\mathcal M_0(\cR_G|_F)$. The product and the adjoint in the above formula are computed in the algebra  $\mathcal M_0(\cR_G|_F)$.

%The  multiplicativity property (\ref{mult})  is proved in \cite{Ra1}.

\end{rem}

\section{A canonical presentation by generators and relations of the algebra $\mathcal M_0(\cR_G|_F)$}

We first introduce the following definition:

 \begin{defn}\label{def2}
Let $\B$ be the minimal $\sigma$-algebra of measurable subsets of $F$ that is left invariant under the action of the partial isomorphisms  considered in formula  (\ref{choice}) and that contains    all the bijectivity domains for the above transformations, which are the sets considered in Lemma \ref{bijective}:
$$\alpha_i^{-1}[\Gamma_g] gF\cap F \quad \mbox{\rm and}\quad g^{-1}[\Gamma_g]\alpha_i F\cap F, g\in G.$$
In the above formula, for each $g\in G$,  the set $(\alpha_i)$ consist of the  right coset representatives for $\Gamma_g$ in $\Gamma$ (i.e.   $\Gamma=\bigcup_i \Gamma_g \alpha_i.$)
\end{defn}

The aim of this section is to provide a description of the $\sigma$-algebra $\B$.
Obviously, with the notation introduced in Definition \ref{def2}, the sets $$g^{-1}[\Gamma_g]\alpha_i F\cap F$$ are equal to the sets
\begin{equation}\label{old1}
\{s\mid gs\in \Gamma_g\alpha_i F\cap F\},
\end{equation}
while the sets
$\alpha_i^{-1}[\Gamma_g] gF\cap F=\alpha_i^{-1}g[\Gamma_{g^{-1}}]F\cap F$ are equal to the sets
\begin{equation}\label{old2}
\{s\in F\mid \alpha_i^{-1}gs\in \Gamma_{g^{-1}}F\} .
\end{equation}
We  decompose  $\Gamma_g$ as a disjoint union of cosets, with respect to smaller normal subgroup
$\Gamma_0\subseteq\Gamma_g$. Thus    $$\Gamma_g=\bigcup\limits_a \gamma_a\Gamma_0.$$
Then the sets in $(\ref{old1})$ and $(\ref{old2})$  are  disjoint unions of sets of the form
$$
\{s\in F\mid gs\in \Gamma_0 F\},\quad g \in G.
$$

We describe explicitly, using the generators constructed  in Definition  \ref{formula}, the action of the partial  isomorphisms, determining   the  equivalence relation  $\mathcal R_G|_F$ on the $\sigma$-algebra $\B$. We will also give an abstract description of the $\sigma$-algebra $\B$.

Let $\mathcal L$ be the the lattice of subgroups of $\Gamma$ generated
by all subgroups of the form $\Gamma_\sigma$, $\sigma \in G$.  We assume that the lattice $\mathcal L$ has a cofinal subset of    normal subgroups in $\Gamma$.

Let
\begin{equation}\label{S}
\mathcal S=\{g\Gamma_0 | \ \Gamma_0\in \mathcal L, g\in G\},
\end{equation}
be the set consisting of all cosets in $G$, corresponding to subgroups in $\mathcal L$. Then $\mathcal S$ has a natural partial relation corresponding to the inclusion relation.
We let
$\mathcal Y$ be the disjoint union
$$\mathcal Y= \bigcup_{n\geq 0}\mathcal S^{\times n}.$$
Here we use the convention that $S^{\times 0}$ is the singleton $\{e\}$, where $e$ is the neutral element of the group $G$.

Let $n$ be a natural number, $g_i$ be arbitrary elements of $G$, $\Gamma_i$ be arbitrary subgroups in $\mathcal L$, and $C_i=g_i\Gamma_i$ for $i=1,2,\ldots,n$.
Let
$$\varepsilon =(C_1,C_2,\dots,C_n)\in \mathcal Y.$$
We let
  \begin{equation}\label{DE}
  D_{\varepsilon}=F\cap C_1F\cap C_2F\cap\ldots\cap C_nF .
 \end{equation}
\noindent  If $n=0$, by convention we let $D_{\{e\}}$  be equal to the set $F$.

Recall that the set $F$ is a measurable fundamental domain, with measure 1, corresponding to  the action of $\Gamma$ on $\cX$. Let  $\Gamma_0$  be a subgroup of $\Gamma_1\cap\ldots \cap\Gamma_n$. It is the obvious that $D_\varepsilon$ is a disjoint union of sets as in equation (\ref{DE}), where
 instead of the subgroups $\Gamma_1,\dots,\Gamma_n$ we may  use only  the subgroup $\Gamma_0$.

Below we will use  the alternative notation
 \begin{equation}\label{notation}
 {C_1}\times {C_2}\times \cdots \times  {C_n}=: F\cap C_1F\cap C_2F\cap\ldots\cap C_nF.
 \end{equation}

 Using the above notation, it obviously follows that the intersection operation, for sets as in  the right hand side of the formula (\ref{notation}), corresponds to  the concatenation operation, if using the notation introduced in the left hand side of the above  formula.
%Let $\overline{C_i}$ be closures, in $\mathcal S$, of the cosets $C_i=g_i \Gamma_i, i=1,2,...,n$.
%
We define $$\nu_F( {C_1}\times {C_2}\times \cdots \times {C_n})= \nu(F\cap C_1F\cap C_2F\cap\ldots\cap C_nF).$$

In the next statement we prove the equality between the $\sigma$-algebra $\B$ and the set consisting of countable unions of sets of the form
$${C_1}\times {C_2}\times \cdots \times {C_n},$$
as in formula (\ref{notation}).

 We will denote by $L^\infty(\B, \nu_F)$ the algebra of bounded $\B$-measurable functions on $\cX$, endowed with measure $\nu_F$ introduced above. By construction  we have
$$L^\infty(\B, \nu_F)\subseteq L^\infty(\cX, \nu).$$

 % Let $ K$ be the profinite completion of the group  $\Gamma$ with respect to the lattice $\mathcal L$ generated by the subgroups of the form
%$\sigma \Gamma \sigma^{-1}\cap\Gamma$, $\sigma \in G$.  Let $\mathcal S$ be the Schlichting completion (\cite{Sch}, \cite{Tz}, \cite{KLM}) of $G$ corresponding to the lattice $\mathcal L$. Then $\mathcal S$ is the disjoint union of the cosets $KgK$, where
%$\Gamma g\Gamma$ runs overs a system of representatives for cosets of $\Gamma$ in $G$. Then $\mathcal S$ is locally compact, totally disconnected group with a canonical Haar measure.

The transformations   $\dot{\wideparen{\Gamma g}}$, generating the equivalence relation $\cR_G|_F$, are subject to the composition equations proved in Theorem
\ref {formula}. This result and the new notation introduced in (\ref{notation}) provide the formulae
\begin{equation}\label{formula1}
      \dot{\wideparen{\Gamma g_1}}\circ \dot{\wideparen{\Gamma g_2}}=
\sum_j \;\dot{\wideparen{\Gamma g_1r_j g_2}}  \chi_{[{(r_jg_2)^{-1}[\Gamma_{g_1^{-1}}}]]},\quad g_1,g_2\in G,
\end{equation}
where $r_j$, $j=1,2,\ldots,[\Gamma:\Gamma_{g_1^{-1}}]$
are right coset representatives for $\Gamma_{g_1^{-1}}$ in $\Gamma$.

We reformulate the content of Lemma \ref{bijective} using the notation from formula (\ref {notation}). Recall that the partial isomorphisms of the set $F$ introduced in the list in formula (\ref{choice}) are identified with partial isometries in the $\ast$-algebra $\mathcal M_0(\cR_G|_F)$ associated to the relation $\cR_G|_F$.

 Fix $g\in G$ and let $(\alpha_i)_i$, $i=1,\dots,[\Gamma :\Gamma_g]$  be a set of right   coset representatives corresponding to  the inclusion $\Gamma_g\subseteq \Gamma$. Then, for all indices $i$,  the  elements
\begin{equation}\label{isometries}
   (\dot{\wideparen{\Gamma g}} )\chi_{[{g^{-1}\Gamma_g\alpha_i }]}\in  \mathcal M_0(\cR_G|_F),
 \end{equation}
     are partial isometries with initial spaces equal to  the projections
$$\chi_{[{g^{-1}\Gamma_g\alpha_i }]}\in L^{\infty} (\B,\nu_F).$$
The range of the isometries in formula (\ref{isometries}) is  equal to the projections $$\chi_{[{\alpha_i^{-1}g\Gamma_{g^{-1}}}]}\in L^{\infty} (\B,\nu_F).$$
The partial inverse (adjoint) of the above isometry is  $\dot{\wideparen{\Gamma g^{-1}\alpha_i}} \chi_{[{\alpha_i^{-1}[\Gamma_g] g}]}$.

%%%%%%%%%%%%%%

\begin{thm}\label{T} Let $\mathcal M_{\rm min}(\cR_G|_F)$
be the $\ast$-subalgebra of the $\ast$-algebra $\mathcal M_0(\cR_G|_F)$,
 % associated to the measurable equivalence  relation $\cR_G|_F$,
  generated by the partial isometries (introduced in formula (\ref{isometries}))
  %  (as in formula (\ref{gens}))
$$(\dot{\wideparen{\Gamma g}} )\chi_{[{g^{-1}\Gamma_g\alpha_i }]},\quad i=1,\ldots,[\Gamma :\Gamma_g],\
g \in G,$$
\noindent  and by the algebra $L^\infty (\B,\nu_F)\subseteq L^\infty(\cX,\nu)$.
Then %The transformations $(\dot{\wideparen{\Gamma g}} )\chi_{[{g^{-1}\Gamma_g\alpha_i }]}$  are measure-preserving partial isomorphisms.
\item (i) The range  of the representation, constructed in Section \ref {outline}, of the Hecke algebra $\mathcal H_0(\Gamma,G) $ into the algebra
$\mathcal M_0(\cR_G|_F)\subseteq B(L^2(F,\nu))$, is contained into the  $\ast$-algebra $\mathcal M_{\rm min}(\cR_G|_F)$.
 \item (ii)
 The action of the  transformations  listed  in formula (\ref{isometries}) on the sets introduced in formula  formula (\ref{notation})
 is  abstractly described as follows:

 Consider an arbitrary element in the set $\Y$:   $${g_1^{-1} \Gamma_0}\times {g_2^{-1}\Gamma_0}\times\cdots\times {g_n^{-1}\Gamma_0}.$$
 % is a finite product in $\mathcal Y$, as in formula (\ref {notation}).
  Let $\Gamma_1$  be a  small enough  (to be determined in the proof) normal subgroup of $\Gamma_g$,  and let $(r_j)_{j=1,\dots [\Gamma_g:\Gamma_1]}$ be a system of right
  coset representatives of the subgroup $\Gamma_1$ in the group $\Gamma_g$.
 Let $(\alpha_i)_{ i=1,\dots, [\Gamma:\Gamma_g]}$ be a system of right cosets representatives  for the subgroup $\Gamma_g$ of $G$.  For each $i$, consider  the corresponding  disjoint  sets decomposition
$$F\cap {g^{-1}\Gamma_g\alpha_i }F=
\bigcup_j [F\cap g^{-1}r_j \alpha_i\Gamma_1F]= \bigcup_j  g^{-1}r_j \alpha_i\Gamma_1.$$

 Then, for all $i,j$, the partial isomorphism  $(\dot{\wideparen{\Gamma g}} )\chi_{[{g^{-1}[\Gamma_g]\alpha_i }]}$  maps the characteristic function of the set
$${g^{-1}r_j \alpha_i\Gamma_1}\times {g_1^{-1} \Gamma_0}\times {g_2^{-1}\Gamma_0}\times\cdots\times {g_n^{-1}\Gamma_0}$$
 onto the characteristic function of the set $${\alpha_i^{-1}r_j^{-1}\Gamma_1 g}
\times  \times_{i=1}^n
{\alpha_i^{-1}r_j^{-1}gg_i^{-1}\Gamma_0}.
 $$

\item (iii). The $\sigma$-algebra $\B$ consists of countable unions of  sets as in formula (\ref{notation}).

\end{thm}

\begin{proof}
Part (i) is a direct consequence of the formula (\ref{heckeformula1}) for the Hecke operators.

To prove part (ii), we note first that
it is sufficient to describe the action of the generators $\dot{\wideparen{\Gamma g}}$, $g\in G$ on a  domain of bijectivity.
For a fixed  $g\in G$, we let  $(\alpha_i)_i$ be  a system
of representatives for $\Gamma_g \subseteq \Gamma$. Thus $\Gamma = \bigcup_i \Gamma_g
\alpha_i$. Consider the restriction of $\dot{\wideparen{\Gamma g}}$
to the bijectivity domain $$g^{-1} \Gamma_g \alpha_i F \cap F = \{s\in F \mid gs\in\Gamma_g \alpha_i F\}.$$

If $\Gamma_0$ is  a normal  subgroup of $\Gamma$, we let   $$A_{g_1,\ldots,g_n,\Gamma_0} = \{s\in F \mid g_i s \in \Gamma_0 F\},$$
be one of the generators of the $\sigma$-algebra $\B$ introduced above.

%%%%%%%%

Recall that  the $\sigma$-algebra $\B$ is generated by subsets of  $F$  of the
form $$A_{g_1,g_2,\ldots,g_n,\,\Gamma_1,\ldots,\Gamma_n}
=
\{s\in F\mid g_1s\in \Gamma_1F,\ldots,g_n\, s\in\Gamma_n F\},
$$
where $g_1,g_2,\ldots,g_n\in G$, $\Gamma_i,\ldots,\Gamma_n$
are subgroups of $\Gamma$ of finite index, in the directed subset $\mathcal L$
of subgroups of $\Gamma$.

It is clear that by dividing each of the subgroups $\Gamma_i$ into cosets with respect to a smaller
common subgroup $\Gamma_0$, we arrive at the situation where we  work with the
$\sigma$-algebra generated by  subsets  of $F$
 of the form
$$
A_{g_1,\ldots,g_n,\Gamma_0} = \{s\in F \mid g_i s \in \Gamma_0F\}.
$$

We may also  reduce to the case in which we work only with  $g_i$ in a fixed system $\mathcal R$  of representatives for
cosets of $\Gamma$ in   $G$.  Let  $r_j$ be a system of representatives
for $\Gamma_0$ in $\Gamma$, thus
$\Gamma = \bigcup_j r_j \Gamma_0$, $j = 1,2,\ldots,[\Gamma :\Gamma_0]$.

Then the following sets:
$$
A_{g_1,\ldots,g_n,r_{j_1},\ldots,r_{j_n},\Gamma_0} = \{s\in F \mid g_i s \in r_{j_i}
\Gamma_0 F,\ i = 1,2,\ldots,n\},
$$
where $g_1,\ldots,g_n$ run over the system of representatives $\mathcal R$ of $\Gamma$ in $G$, $j_1,\ldots,j_n$ run over $\{1,2,\ldots,[\Gamma:\Gamma_0]\}$ and $\Gamma_0$ runs over $\mathcal L$, and $n\in\mathbb N$,  are a system of generators of the $\sigma$-algebra $\B$ of subsets of $F$ .

With $\alpha_i$ as in the statement of the theorem, $\dot{\wideparen{\Gamma g}}$ is bijective on the set
\begin{equation}\label{old3}
\{s\in F \mid gs \in \Gamma_g \alpha_i F\}.
\end{equation}
\noindent
We let $\Gamma_1$  be a small enough normal subgroup. The conditions on how small the group $\Gamma_1$ has to be taken
will be determined later in the proof. We decompose
$\Gamma_g$ as a disjoint union $$\bigcup_{j=1}^{[\Gamma_g:\Gamma_1]} r_j \Gamma_1,$$ of left   cosets of the subgroup $\Gamma_1$.

Hence the  set in formula (\ref{old3}) becomes the disjoint union of the sets
$$
\{s\in F \mid gs \in r_j \alpha_i\Gamma_1 F\}.
$$

We want to  determine the image through $\dot{\wideparen{\Gamma g}}$ of the set
\begin{equation}\label{old4}
\{s\in F \mid gs \in r_j \alpha_i\Gamma_1 F\} \cap
\{s\mid g_i s \in \Gamma_0 F\}.
\end{equation}
Note that the normality of $\Gamma_1$ yields
$$r_j\alpha_i\Gamma_1=\Gamma_1r_j\alpha_i.$$

We fix an element $f$ in the set (\ref{old4}).  Then $gf$ is of the form $\theta_1r_j\alpha_if_1$
with $\theta_1\in\Gamma_1$ and $f_1\in F$. Moreover, $$g_if\in\Gamma_0 F.$$
Then
$$\dot{\wideparen{\Gamma g}}(f)=f_1,$$
and
$$f=g^{-1}\theta_1 r_j\alpha_i f_1.$$
The condition that $g_i f\in\Gamma_0 F$ then translates into
$g_i(g^{-1}\theta_1 r_j\alpha_i) f_1\in\Gamma_0 F$,
which is the same as
$$
f_1\in\alpha_i^{-1}r_j^{-1}\theta_1^{-1}gg_i^{-1}\Gamma_0F.
$$
We choose the subgroup  $\Gamma_1$ to be   small enough so that there exists a fixed
subgroup $\Gamma_2$ of $\Gamma_0$ such that, for all $i$,
$$\Gamma_1gg_i^{-1}=gg_i^{-1}\Gamma_2 .
$$
Hence the condition on $f_1$ is that
$$
f_1\in\alpha_i^{-1}r_j^{-1}gg_i^{-1}\Gamma_0F.
$$
We also have to impose the condition that $f_1$ belongs to the image of
$$\{s\mid gs\in r_j\alpha_i\Gamma_1 F\}$$ through $\dot{\wideparen{\Gamma g}}$.
But for all $j$, we have $$f=g^{-1}\theta r_j\alpha_if_1.$$
Hence
$$f_1=\alpha_i^{-1}r_j^{-1}\theta^{-1}gf,$$
and hence
$$f_1\in\alpha_i^{-1}r_j^{-1}\Gamma_1 gF\cap F\subseteq\alpha_i^{-1}\Gamma_g g\Gamma\cap F.$$
Thus the necessary and sufficient  condition on $f_1$ is that it belongs to
$$
\{s\in F\mid r_j\alpha_i s\in \Gamma_1 g F\}\cap\bigcap_{i=1}^n
\{g_ig^{-1}r_j\alpha_i s\in\Gamma_0 F\}.
$$
Note that the first set in the intersection  is also described by the formula $$\{s\in F\mid g^{-1}r_j\alpha_i s\in g\Gamma_1g^{-1} F\},$$
where $g\Gamma_1g^{-1}\subseteq \Gamma$ since $\Gamma_1\subseteq\Gamma_{g^{-1}}$.

When translating this in terms of the measure space $L^{\infty}(\B, \nu_F)$, we obtain part (ii) in the theorem.

Because of (ii), the $\sigma$-algebra generated by measurable  subsets of $F$,  as considered  in formula (\ref{notation}),  is left invariant by the transformations listed in formula (\ref{isometries}).
Hence, by the Definition \ref{def2} this algebra coincides with the $\sigma$-algebra $\B$.
This concludes the proof of part (iii) of the theorem.
\end{proof}

\section{The situation $G=\PGL_2(\Z[\frac1p])$, $\Gamma=\PSL_2(\Z)$, $p$ odd prime}

In this section  we assume that   $G=\PGL_2(\Z[\frac1p])$ and $\Gamma=\PSL_2(\Z)$, where $p\geq 3$ is a prime number. Recall that $\Z[\frac1p]$ is the ring of rational numbers whose denominators are  powers of the prime number $p$.
We  use the properties of the generators of the algebra $\mathcal M_0(\cR_G|_F) $ associated to the equivalence relation $\cR_G|_F$, and their multiplication formula, proved in  Theorem \ref {formula}.

 We  prove  that in this case  the measurable, countable, ergodic equivalence relation ${\cR_G|_F}$ on the fundamental domain $F$  is  induced by an a.e. free, ergodic,  probability, measure-preserving action of the free group on $(p+1)/2$ generators on $F$.  The partial isomorphisms of the set $F$,  introduced in Definition \ref{pi} and Lemma  \ref{generators}, provide an explicit L-graphing (\cite{Ga}, see also  \cite{Ke1}, Section 17, page  58) of minimal cost $(p+1)/2$, for the  treeable equivalence relation  ${\cR_G|_F}$.

% The in such a way that the Hecke algebra corresponding to the inclusion $\Gamma \subseteq G$ of Hecke operators acting on $L^2(F, \nu)$ coincides with the image, through the Koopmann representation of the free group, of the radial algebra in the free group.

 We are indebted to the anonymous referee, of a first version of this paper, for pointing out to the author that a related argument was considered in the paper \cite{Ad}.

In the case considered in this section,
the structure of the Hecke $\ast$-algebra  $\mathcal H_0(G,\Gamma)$  is
well known. We recall from \cite{Kr} the basic facts about the structure of this algebra and its action on the linear space $\ell^2(\Gamma\backslash G)$ freely generated by the set $\Gamma\backslash G$ of cosets of $\Gamma$ in $G$.

For $n\in \mathbb N$, consider the following group element in $G$:
$$\sigma_{p^n}=\left(
\begin{matrix}
p^n&0\\0&1
\end{matrix}\right).$$

The double cosets of $\Gamma$ in $G$ are
$$\{[\Gamma \sigma_{p^n}\Gamma] \mid n\in \mathbb N\}.$$
The identity element of the Hecke algebra is $[\Gamma]$, the coset corresponding to the identity element in $G$.
The Hecke algebra is commutative and the multiplication rule in the Hecke algebra is  as follows:
\begin{equation}\label{multrad}
[\Gamma \sigma_{p^n}\Gamma] [\Gamma \sigma_{p}\Gamma]=
[\Gamma \sigma_{p^{n+1}}\Gamma]+p [\Gamma \sigma_{p^{n-1}}\Gamma],\quad  n\geq 1,
\end{equation}
and
\begin{equation}\label{multrad1}
[\Gamma \sigma_{p}\Gamma] [\Gamma \sigma_{p}\Gamma]=
[\Gamma \sigma_{p^2}\Gamma] +(p+1) [\Gamma].
\end{equation}

The $\ast$-operation on the Hecke algebra is
\begin{equation}\label{selfa}
\big([\Gamma \sigma_{p^n}\Gamma]\big)^\ast= [\Gamma (\sigma_{p^n})^{-1}\Gamma]= [\Gamma \sigma_{p^n}\Gamma].
\end{equation}
There are exactly $p+1$ (left or right) cosets of  $\Gamma$ contained in the double coset $[\Gamma \sigma_{p}\Gamma]$, which we enumerate as follows:
\begin{equation}\label{enumeration}
\Gamma \sigma_{p}\Gamma=\bigcup_{i=1}^{p+1}\Gamma\sigma_p s_i.\end{equation}

The above sets of formulae (\ref{multrad}), (\ref{multrad1}), (\ref{selfa}) is identical to the relations that hold  in the radial algebra of a free group. We briefly recall below the structure of the radial algebra in the group algebra $\mathbb C(F_{(p+1)/2})$ of a free
group with ${(p+1)/2}$ generators (see e.g. \cite{PS}). For $n\in \mathbb N$, let $\chi_n\in \mathbb C(F_{(p+1)/2})$ be the sum of all words of length
$n$. In particular $\chi_0$ is the identity element, and the selfadjoint
elements $\chi_n$ verify the equalities
(\ref{multrad}) and  (\ref{multrad1}), with $\chi_n$ substituting $[\Gamma \sigma_{p^n}\Gamma]$ in the above formulae, for all $n\in\mathbb N$.
The $C^\ast$-algebra generated by $\chi_n$, $n\in {\mathbb N}$ in the reduced $C^\ast$-algebra of the free group is called the radial algebra.

The formula for the Hecke operator
$$T_p=:T([\Gamma \sigma_{p}\Gamma]),$$
corresponding to the double coset $[\Gamma \sigma_{p}\Gamma]$ and   acting on $L^2(F,\nu)$, associated to the Koopman unitary representation $\piK$ of $G$ into the unitary group of
\ $L^2(\cX,\nu)$ (see the introductory section for definitions), is
\begin{equation}\label{hecke1}
(T_p )f(x)= \sum_{i=1}^{p}
f((\sigma_p s_i)^{-1}x),\quad x\in \cX,  f\in L^2(F,\nu).
\end{equation}
The Hecke operators $T_{p^n}=:T([{\Gamma\sigma_{p^n}\Gamma}])$, for $n>1$ are described by a similar formula.

We also recall bellow a few facts (\cite{Ke1}) about special sets of generators for the countable, measurable equivalence relation  ${\cR_G|_F}$. Consider a set
$$\Phi_0=\{\phi \mid \phi :A_\phi\rightarrow B_\phi\},$$
\noindent consisting of measure preserving, partial isomorphisms $\phi$ of $F$, with domain $A_\phi$ and range $B_\phi$. Note that it automatically follows that the sets  $A_\phi$ and
$B_\phi$ have equal measure.

The set $\Phi_0$ is called a L-graphing for the measurable  equivalence relation ${\cR_G|_F}$ if,  with the exception of sets of measure 0,
two points $x,y$ in $F$ are equivalent with respect to ${\cR_G|_F}$ if and only if there exist
$\phi_1,\phi_2,\dots, \phi_n$ in $\Phi_0$,  and $\varepsilon_1, \varepsilon_2, \dots, \varepsilon_n$ in $\{\pm 1\}$ such that
$$y=\phi_1^{\varepsilon_1}\circ \phi_2^{\varepsilon_2} \circ \dots \circ \phi_n^{\varepsilon_n}(x).$$

An L-graphing $\Phi_0$ is called an L-treeing (see \cite{Le}, \cite {Ga} or see \cite {Ke1}, Section 21, page 62), which is also  simply referred to as to  a treeing for the relation ${\cR_G|_F}$, if the choice is unique in a relation as above.

More precisely, an L-graphing $\Phi_0$ is a treeing, if almost everywhere,
an equality  the form
$$x=\phi_1^{\varepsilon_1}\circ \phi_2^{\varepsilon_2} \circ \dots \circ \phi_n^{\varepsilon_n}(x),$$ that holds true on a measurable set, not of measure 0, implies that
the word
$$\phi_1^{\varepsilon_1}\circ \phi_2^{\varepsilon_2} \circ \dots \circ \phi_n^{\varepsilon_n},$$
is trivial, due to cancellations of the elements in $\Phi_0$ with their inverses.

The cost $C(\Phi_0)$ (\cite{Le}) of an L-graphing $\Phi_0$  is equal to the sum of the measures of the domains of the partial isomorphisms in the set $\Phi_0$. If the L-graphing $\Phi_0$ is a treeing, then the equivalence  relation is treeable (for definitions see \cite{Le},\cite {Ad}, \cite {Ga}). In this  latest case, the cost $C(\Phi_0)$ is equal to the cost $C(\cR_G|_F)$ of the relation
$\cR_G|_F$ (\cite{Ga}). If the relation $\cR_G|_F$ is treeable of integer cost $n$, then there exists an ergodic, measure preserving action of a free group $F_n$ on $F$, whose orbits are exactly the equivalence classes of the countable, measurable equivalence relation $\cR_G|_F$
(\cite {Hj}).

\begin{thm} \label{prime}
%
%Thus, two points in $F$ are equivalent, with respect to countable equivalence relation $\cR_G|_F$, if and only if
%the two points belong to   the same orbit of $G$ acting on $\cX$.
%
Let $\Gamma g$ be any of the cosets in the decomposition of the double coset $\Gamma\sigma_p\Gamma$ in left cosets,  as in  formula (\ref{enumeration}). For each such $\Gamma g$, let $(s_i^g)_{i=1}^{p+1}\subseteq \Gamma $  be a set of right coset representatives for the subgroup $\Gamma_g$ of $\Gamma$.

%$$\Gamma=\bigcup_{i=1}^{p+1}\Gamma_gs_i^g,$$ where the elements
%$(s_i^g)_{i=1}^{p+1}$ are the coset representatives.

Consider  the subset
\begin{equation}\label{gens1}
\mathcal A_p=\{\dot{\wideparen{\Gamma g}}|_{[g^{-1}[\Gamma_g]s_i^gF\cap F]}\mid \Gamma g \subseteq
%{\rm \ runs \ through\   the\ cosets\ in\ }
 \Gamma\sigma_p\Gamma, i=1,2,\dots, p+1
 % {\rm \ and\ }  \Gamma=\bigcup[\Gamma_g]s_i^g
 \},
 \end{equation}
of the set $\Phi_{\cR_G|_F}$ of partial isomorphisms introduced in Lemma \ref {generators}.

Then
\item (i)
The  set $\mathcal A_p$  is closed with respect to  the inverse operation for partial isomorphisms. No element in $\mathcal A_p$ is equal to its inverse.
 %This holds true, in general, by simply imposing (\cite{Kr}) that
%$\Gamma\sigma\Gamma=\Gamma\sigma^{-1}\Gamma$, for all $\sigma \in G$.
\item (ii) Let $\mathcal  {A}^0_p\subseteq  \mathcal A_p$ be a subset of representatives (modulo the inverse operation) for the partial isomorphisms in the set $\mathcal A_p$. This choice  means that if  $\phi \in \mathcal  {A}^0_p$ then  $\phi \notin \mathcal  {A}^0_p$.

 Then, the set $\mathcal A_p$ is a canonical   L-graphing for the equivalence relation $\mathcal R_G|_F$.  The set $\mathcal  {A}^0_p$ is a treeing for the measurable equivalence relation
  $\mathcal R_G|_F$.
 %This means that two points $x, y$ are equivalent in $\mathcal R_G|_F$ if and only if there exists  a chain
 %of elements $\dot{\wideparen{\Gamma g_1}},\dot{\wideparen{\Gamma g_2}}\ldots,\dot{\wideparen{\Gamma g_n}}$, restricted to their injectivity domains as in Lemma \ref {bijective}, belonging   to $\mathcal A_p$, such that
% \begin{equation}\label{graphing}
% \dot{\wideparen{\Gamma g_1}}\circ\dot{\wideparen{\Gamma g_2}}\circ\cdots\circ \dot{\wideparen{\Gamma g_n}}\ x=y
% \end{equation}
%Note that the set $\mathcal A$ has cardinality $p+1$.
%Moreover, any equality of the form
%$$\dot{\wideparen{\Gamma g_1}}\circ\dot{\wideparen{\Gamma g_2}}\circ\cdots\circ\dot{\wideparen{\Gamma g_n}}\ f=f,\quad  f\in F,$$
%where the transformations $\dot{\wideparen{\Gamma g_i}}$ are restricted to an injectivity domain as in Lemma \ref{bijective}, holds true, if and only the composition  of the transformations $\dot{\wideparen{\Gamma g_i}}$,  in the above formula,
%is the identity transformation, restricted to the injectivity domain of the last transformation in the chain.
%Consequently, the graphing  $\mathcal A_p$, is a treeing (\cite{Ga}) for the equivalence
%relation $\cR_G|_F$.
%The total measure of the domains of the elements
%in $\mathcal A_p$ is $(p+1)$.
 \item (iii)
 In particular, the equivalence relation
$\cR_G|_F$ is treeable (see  \cite{Ad}), and the  cost (\cite{Ga}) of the relation $\cR_G|_F$ is $(p+1)/2$.
\end{thm}

%In addition, we can arrange that that generators of $F_{\frac{p+1}{2}}$ are built
%from pieces of the transformations of $\dot{\wideparen{\Gamma g}}$, $\Gamma g \subseteq \Gamma \sigma_p \Gamma$, glued together into bijective transformations.
%Hence the radial elements in group algebra of the group $F_{\frac{p+1}{2}}$ (that is the selfadjoint elements $\chi_n$ = sum of words in
%the generators of length $n$, $n\in \N$) have the property that $\chi_n$ as an operator
%on $L^2(F)$ coincides with
 %the Hecke ope\-rator~$T_{\sigma_{p^n}}$, where $$\sigma_{p^n}=\left(
%\begin{matrix}
%p^n&0\\0&1
%\end{matrix}\right).$$

\begin{proof}
The fact that the set $\mathcal A_p$ is closed under the inverse operation is a consequence of the property expressed in formula (\ref{selfa}) and of the construction in Lemma \ref {generators}. This proves part (i).

We next prove part (ii).   Recall (\cite{Serre}) that the action of $\Gamma\sigma_p \Gamma$ on the
cosets in $\Gamma \setminus G_p$ copies exactly the action of the radial algebra on the
elements of the free group $F_{(p+1)/2}$.
Consequently, the Cayley tree of $F_{(p+1)/2}$ with origin at the neutral element $e$ is identified with the set of cosets
$\Gamma \setminus G$. The elements at distance $n$ to the origin correspond to the cosets in
$\Gamma\sigma_{p^n} \Gamma$, $n\in\mathbb N$. Let
$$\chi_1 = \sum\limits_{i=1}^{(p+1)/2} s_i + s_i^{-1}, $$
where $(s_i)_{i=1}^{(p+1)/2}$ are the generators of $F_{(p+1)/2}$. Then $\chi_1$ is the radial element of order 1 (the discrete laplacian) on $F_{(p+1)/2}$. The action of $\chi_1$ on $F_{(p+1)/2}$ corresponds bijectively to multiplication by
$\Gamma\sigma_p \Gamma$ in the space of cosets. More precisely (\cite{Serre}), there exists a bijection
$$\Psi: F_{(p+1)/2}\to \Gamma \setminus G,$$
such that $\Psi$ maps the set of words of length $n$ in $F_{(p+1)/2}$ into the set of left cosets of $\Gamma$ contained in  the double coset $\Gamma \sigma_{p^n}\Gamma$, for all $n\in \mathbb N$.
Moreover, for every $w\in  F_{(p+1)/2}$, the
 set $$\Psi (\{s_iw,\ s_i^{-1}w,\mid i=1,2,\ldots,(p+1)/2\}),$$
is equal to   the set of cosets in $[\Gamma\sigma_p \Gamma]\Psi(w)$.
% for any word $w$ in  $F_{\frac{p+1}{2}}$.

The fact that $\mathcal A_p$ is a L-graphing   follows from
 Lemma  \ref {generators} and from the fact that, because of formula (\ref{multrad})  in the introductory part of the introduction, the coset $\Gamma \sigma_p\Gamma$ is a selfadjoint  generator for the Hecke algebra $\mathcal H_0(\Gamma, G)$.

Let $\Gamma g_1,\Gamma g_2,\ldots,$ $\Gamma g_n$ be  cosets in $[\Gamma\sigma_p \Gamma]$. Assume that the equality
\begin{equation}\label{product1}
\dot{\wideparen{\Gamma g_1}}\circ\dot{\wideparen{\Gamma g_2}}\circ\cdots\circ\dot{\wideparen{\Gamma g_n}}\  f=f
\end{equation}
holds true for $f$ in a measurable subset of non-zero measure of $F$.

In the above equation, we assumed that all the  transformations $$\dot{\wideparen{\Gamma g_1}},\ldots,\dot{\wideparen{\Gamma g_n}},$$ are restricted to a domain of injectivity as in Lemma \ref{bijective}.
By Theorem \ref{formula}, this corresponds   to an equality of the form,
$$
\gamma g_1 r_1 g_2 r_j \ldots g_{n-1} r_{n-1} g_n f = f,
$$
where $\gamma, r_1,r_2,\dots, r_{n-1}$ all belong to $\Gamma$.

Because of the properties of the map  $\Psi$ constructed above, an equality as above  may hold true for $f$ in  a non-zero measure subset of $F$ if and only if we can identify   successive cancellations of a transformation in the set $\mathcal A_p$ by its inverse. 

Indeed,  these    cancellations correspond, because of the properties of the map $\Psi$,  to products      $g_{j_1} r' g_{j_2}$, in the above formula, that belong to the group $\Gamma$.
Thus an identity as    in formula (\ref{product1}) corresponds to the fact that we are successively   multiplying the transformations $\dot{\wideparen{\Gamma g_j}}$ in the product from the formula (\ref{product1}), restricted to an injectivity domain as in Lemma \ref{bijective}, with their inverse transformation, until we obtain the identity element.

Hence the product of the the transformations in formula (\ref{product1}) is the identity (restricted to an injectivity domain corresponding to the first transformation in the product).

Thus the equivalence relation is treeable, with an L-treeing   consisting in a set of representatives, modulo the inverse operation, of the partial isomorphisms in the  set $\mathcal A_p$.

Since the set $\mathcal A_p$  is closed under the inverse operation and the sum of the  measures  of the domains of the partial isomorphisms in $\mathcal A_p$  is $p+1$, the cost of the relation must be $\frac{p+1}{2}$.
\end{proof}
%By Hjorth theorem \cite{Hj}, we can find a free group $F_{\frac{p+1}{2}}$ whose orbits define
%the relation $\cR_{G_p}|F$.

% Since we have that the set of  point images and set of  point preimages of  the $p+1$ elements   elements in the
%set $\mathcal A$ have cardinality exactly $p+1$ in the set $F$, and $\mathcal A$ is closed to the inverse operation, it follows that there is no obstruction into measurably breaking the elements in $\mathcal A$ into partial transformations, and recomposing them into a set $\mathcal A'$ of $p+1$ bijective transformations, such that the set
%$\mathcal A'$ is closed to the operation of taking the inverse. The new %elements in the set $\mathcal A'$ preserve the property of having no relations
%(except the cancellations due to the inverse operation), and hence the elements in $\mathcal A'$ generate a free group action. The radial algebra element $\chi_1$ is the sum of the transformation operators induced by the elements in  the set  $\mathcal A'$ and this is equal to the  corresponding sum for the elements in $\mathcal A$.
%But this sum  is the Hecke operator corresponding to the double coset $[\Gamma \sigma_p \Gamma]$, acting on $F$.

\begin{rem}
We recall from above  that by Hjorth's theorem (\cite{Hj}), the previous result proves that there exists a free group factor $F_{(p+1)/2}$ acting
freely on $F$, whose orbits are the equivalence classes of the relation  $\cR_G|_F$.

Note the Hecke operator $T_p$ introduced in formula (\ref{hecke1}) has the formula (see formula (\ref{heckeformula}) in Section \ref{product})
$$T_p=\sum _{\phi\in \mathcal A_p} \phi\in \mathcal A_{\operatorname{Koop}}(\cR_G|_F),$$
where $\mathcal A_{\operatorname{Koop}}(\cR_G|_F)$ is the $C^\ast$-algebra introduced
in formula (\ref{akoop}).
%Here, as explained in the introduction, we identify $\phi$ with the partial isometry $U_\phi$, belonging to  the algebra $\mathcal M_0(\cR_G|_F)$, and defined in  formula (\ref{uphi}).

Because of Hjorth's theorem, the algebra $\mathcal A_{\operatorname{Koop}}(\cR_G|_F)$ is isomorphic to the Koopman representation (\cite{Ke}) on $L^2(F,\nu)$ of the crossed product algebra $$C^\ast(F_{(p+1)/2})\rtimes L^\infty(F,\nu)).$$

Consequently, the $C^\ast$-algebra $\mathcal A_{\operatorname{Koop}}(\cR_G|_F)$ contains two algebras which are $\ast$-isomorphic to the radial algebra
of the free group $F_{(p+1)/2}$. One of these algebras is the image, in the Koopman representation, of   the radial algebra of the free group itself, which acts as above on $F$.

 The second algebra is generated by the representation of the Hecke operators in $\mathcal A_{\operatorname{Koop}}(\cR_G|_F)$. The embedding of the Hecke algebra in $\mathcal A_{\operatorname{Koop}}(\cR_G|_F)$ is constructed in the introduction.

We conjecture that the two representations of the radial algebra are unitarily equivalent in $\mathcal A_{\operatorname{Koop}}(\cR_G|_F)$. Moreover, we conjecture that the Koopman unitary representation of the free group
$F_{(p+1)/2}$ is tempered, in the sense considered by Kechris (\cite{Ke}). The second conjecture is suggested by the results in the paper
\cite{WW}. A simultaneous  positive answer to these conjectures would provide a positive answer to the Ramanujan-Petersson conjectures.

Note that because of the above, the Hecke operators associated with $\Gamma \subseteq G$ and a measure-preserving action of $G$ on an infinite measure space, admitting a fundamental domain for the action of $\Gamma$, are defined  similarly to the Hecke operators considered in \cite{LPS}.
\end{rem}

We introduce  the following  definition that is suggested by   the preceding theorem and by the notion of weak orbit equivalence (\cite{Fu}).

\begin{defn}\label{ior}
Let $H_1,H_2$ be two discrete groups. We will say that $H_1$
is  infinitesimally orbit equivalent to $H_2$
if there exist an infinite, ergodic, measure-preserving, a.e. free
action of $H_2$ on an infinite measure space $Y$, and $F$ a finite measure subset of $Y$, such that if
$\cR_{H_2}$ is the countable equivalence relation induced on $Y$ by
the orbits of $H_2$, then $\cR_{H_2}|_F$ is orbit equivalent to the measurable equivalence relation corresponding to a   free, ergodic, measure-preserving action
of $H_1$ on $F$.
\end{defn}

Thus, in particular we proved the following:

\begin{cor}
The group $F_{(p+1)/2}$ is infinitesimally orbit equivalent
to the group  $\PGL_2(\Z[\frac1{p}])$.
\end{cor}

\section{Matrix coefficients of Hecke operators}

In this section we compute the matrix coefficients of the Hecke operators, by a method related to the previous considerations. The context is the same as in Section \ref{product}.

We let $G$ be a countable discrete group
acting ergodicaly, by measure-preserving transformations, on
an infinite measure space $(\cX,\nu)$ with $\sigma$-finite
measure $\nu$. In addition, we assume that $\Gamma \subseteq G$ is an almost
normal subgroup that has a fundamental domain $F$ of finite measure 1 in
$\cX$. Let $\pi=\piK$ be the Koopman (see e.g. \cite{Ke}) representation of $G$ into $L^2(\cX,\nu)$.

On the $\Gamma$-invariant functions on $\cX$ we introduce the (Petersson) scalar product
given by integration over $F$. Consequently, the Hilbert space of \ $\Gamma$-invariant functions on $\cX$, denoted by $L^2(\cX,\nu)^\Gamma$, is identified with the Hilbert space $L^2(F,\nu|_F)$. Recall that  for $\sigma\in G$, we denote by $T(\Gamma \sigma \Gamma)$ the Hecke operator on $L^2(F,\nu)$ associated to the double coset $\Gamma \sigma \Gamma$.

\begin{prop} \label{last}
We use the notation and definitions introduced above.
 We assume that the characteristic function $\chi_F$ is cyclic for $\piK$ and consider the $\Gamma$ -invariant  function (depending on the coset $\Gamma \sigma$) defined by
 $$\tilde{\chi_{\Gamma \sigma}}=\sum_{\gamma \in \Gamma}\pi(\gamma)\chi_{\sigma F}.$$

\noindent
 Let $\mathcal S$ be the linear subspace of $L^2(\cX,\nu)^\Gamma$ generated by the functions
 $\tilde{\chi_{\Gamma \sigma}}$, where $\Gamma \sigma$ runs over all cosets of $\Gamma$ in $G$. By the  assumptions   on the function $\chi_F$ , the space  $\mathcal S$ is dense in $L^2(\cX,\nu)^\Gamma$.

Then we have the following formulae:
\begin{equation}\label{old5}
 T^{\Gamma\sigma \Gamma} \tilde{\chi_{\Gamma \sigma_1}}=\sum_{\Gamma \theta \in [\Gamma\sigma \Gamma][\Gamma \sigma_1]}  \tilde{\chi_{\Gamma \theta}},
 \quad \sigma, \sigma_1 \in G,
\end{equation}
where $[\Gamma \theta ]$ runs over the set of left cosets in the   product  $[\Gamma\sigma \Gamma][\Gamma \sigma_1]$.

The Hilbert space scalar product in $L^2(\cX,\nu)^\Gamma$ is computed by
\begin{equation}\label{old6}
 % \alpha(\Gamma\sigma_1, \Gamma\sigma_2)=
  \langle (\tilde{\chi_{\Gamma \sigma_1}},\tilde{ \chi_{\Gamma \sigma_2}}\rangle_{L^2(\cX,\nu)^\Gamma}=
  %$$
  %$$=
  \sum_{\gamma\in \Gamma}\nu(\sigma_1^{-1}\gamma\sigma_2 F\cap F),\quad  \sigma_1, \sigma_2 \in G.
  \end{equation}
\end{prop}

  Before proving the proposition, we note that formulae (\ref{old6}) may be used to define a direct scalar product on the linear space of cosets 
  ${\mathbb C}(\Gamma \sigma| \sigma \in G)$ freely generated by the set of cosets of $\Gamma$ in $G$. Note that  $\tilde{\chi_{\Gamma \sigma}}$ is not a necessarily a characteristic function, although it is a step function.

Formula (\ref{old5}) proves that the action of the Hecke operators on $L^2(F,\nu)$ copies the action of the Hecke algebra on $L^2(\Gamma \backslash G)$; the only difference is showing up in the formula of the scalar product, which is now given by a new,   positive definite scalar product on the space of cosets.

Obviously, if $\sigma_1$ is the identity element in formula (\ref{old6}), then the value of the scalar product is 1. Also, if we decompose the coset $\sigma_1\Gamma\sigma_2$, $\sigma_1,\sigma_2\in G$ as a finite  disjoint union  $\cup_j g_j\Gamma_{h_j}$ of cosets of smaller  subgroups $\Gamma_{h_j}$ with $g_j, h_j \in G$, then the computation of the right hand side  term in formula (\ref{old5}) is reduced to the calculation of the distribution
$$\nu (\sigma_1F \cap s_i\Gamma_{\sigma_2}F),\quad \sigma_1,\sigma_2 \in G,$$
where $s_i$ are the right coset representatives for the group $\Gamma_{\sigma_2}$ in $\Gamma$.

\begin{proof} [Proof of Proposition \ref{last}]
Formula (\ref{old5}) is also a consequence of the considerations in Appendix 2 in (\cite{Ra}). To prove formula (\ref{old5}), we note the equality
$$\langle (\tilde{\chi_{\Gamma \sigma_1}}, \tilde{\chi_{\Gamma \sigma_2}}\rangle_{L^2(\cX,\nu)^\Gamma}=
\sum_{\gamma_1,\gamma_2\in \Gamma} \nu(F\cap\gamma_1\sigma_1 F\cap \gamma_2\sigma_2F),\quad
\sigma_1,\sigma_2\in G.$$
Since $F$ is a fundamental domain, and $\nu$ is a $G$-invariant measure, this sum is further equal to
$$\sum_{\gamma\in\Gamma}\nu (\sigma_1F\cap\gamma\sigma_2 F).$$
From here we deduce formula (\ref{old6}). If $\sigma_1$ is the identity, the formula computes the measure of $\sigma_2F$, which is 1 by hypothesis. The rest of the statement is obvious.
\end{proof}

{\bf Acknowledgment} The author is deeply indebted to the anonymous referee of a first version of this paper for various comments and for making the author  aware of the reference \cite{Ad}. The author is also grateful to Florin Boca and  to the second anonymous   referee of the second version of this paper for very useful  comments.

\end{document}